# GENEALOGICAL PARTICLE ANALYSIS OF RARE EVENTS

By Pierre Del Moral and Josselin Garnier

*Université de Nice and Université Paris 7*

In this paper an original interacting particle system approach is developed for studying Markov chains in rare event regimes. The proposed particle system is theoretically studied through a genealogical tree interpretation of Feynman–Kac path measures. The algorithmic implementation of the particle system is presented. An estimator for the probability of occurrence of a rare event is proposed and its variance is computed, which allows to compare and to optimize different versions of the algorithm. Applications and numerical implementations are discussed. First, we apply the particle system technique to a toy model (a Gaussian random walk), which permits to illustrate the theoretical predictions. Second, we address a physically relevant problem consisting in the estimation of the outage probability due to polarization-mode dispersion in optical fibers.

**1. Introduction.** The simulation of rare events has become an extensively studied subject in queueing and reliability models [16], in particular in telecommunication systems. The rare events of interest are long waiting times or buffer overflows in queueing systems, and system failure events in reliability models. The issue is usually the estimation of the probability of occurrence of the rare event, and we shall focus mainly on that point. But our method will be shown to be also efficient for the analysis of the cascade of events leading to the rare event, in order to exhibit the typical physical path that the system uses to achieve the rare event.

Standard Monte Carlo (MC) simulations are usually prohibited in these situations because very few (or even zero) simulations will achieve the rare event. The general approach to speeding up such simulations is to accelerate the occurrence of the rare events by using importance sampling (IS) [16, 24]. More refined sampling importance resampling (SIR) and closely related sequential Monte Carlo methods (SMC) can also be found in [4, 10]. In all of these well-known methods the system is simulated using a new set of input probability distributions, and unbiased estimates are recovered by multiplying the simulation output by a likelihood ratio. In SIR and SMC these ratio weights are also interpreted as birth rates. The tricky part of these Monte Carlo strategies is to properly choose the twisted distribution. The user is expected to guess a more or less correct twisted distribution; otherwise these algorithms may completely fail. Our aim is to propose a more elaborate and adaptative scheme that does not require any operation of the user.

Recently intensive calculations with huge numerical codes have been carried out to estimate the probabilities of rare events. We shall present a typical case where the probability of failure of an optical transmission system is estimated. The outputs of these complicated systems result from the interplay of many different random inputs and the users have no idea of the twisted distributions that should be used to favor the rare event. This is in fact one of the main practical issues to identify the typical conjunction of events leading to an accident. Furthermore these systems are so complicated that it is very difficult for the user, if not impossible, to modify the codes in order to twist the input probability distributions. We have developed a method that does not require twisting the input probability distribution. The method consists in simulating an interacting particle system (IPS) with selection and mutation steps. The mutation steps only use the unbiased input probability distributions of the original system.

The interacting particle methodology presented in this paper is also closely related to a class of Monte Carlo acceptance/rejection simulation techniques used in physics and biology. These methods were first designed in the 1950s to estimate particle energy transmission [15], self-avoiding random walks and macromolecule evolutions [23]. The application model areas of these particle methods now have a range going from advanced signal processing, including speech recognition, tracking and filtering, to financial mathematics and telecommunication [10].

The idea is the following one. Consider an $E$-valued Markov chain $(X_p)_{0 \leq p \leq n}$ with nonhomogeneous transition kernels $K_p$. The problem consists in estimating the probability of occurrence $P_A$ of a rare event of the form $\{V(X_n) \in A\}$ where $V$ is some function from $E$ to $\mathbb{R}$. The IPS consists of a set of $N$ particles $(X_p^{(i)})_{1 \leq i \leq N}$ evolving from time $p = 0$ to $p = n$. The initial generation at $p = 0$ is a set of independent copies of $X_0$. The updating from the generation $p$ to the generation $p + 1$ is divided into two stages:

(1) The selection stage consists in choosing randomly and independently $N$ particles amongst $(X_p^{(i)})_{1 \leq i \leq N}$ according to a weighted Boltzmann–Gibbs particle measure, with a weight function that depends on $V$. Thus, particles with low scores are killed, while particles with high scores are multiplied. Note that the total number of particles is kept constant.

(2) The mutation step consists in mutating independently the particles according to the kernel $K_p$. Note that the true transition kernel is applied, in contrast with IS.

The description is rough in that the IPS actually acts on the path level. The mathematical tricky part consists in proposing an estimator of the probability $P_A$ and analyzing its variance. The variance analysis will provide useful information for a proper choice of the weight function of the selection stage.

The analysis of the IPS is carried out in the asymptotic framework $N \gg 1$ where $N$ is the number of particles, while the number $n$ of mutation–selection steps is kept constant. Note that the underlying process can be a Markov chain $(\tilde{X}_p)_{0 \leq p \leq \tilde{n}}$ with a very large number of evolutionary steps $\tilde{n}$. As the variance analysis shows, it can then be more efficient to perform selection steps on a subgrid of the natural time scale of the process $\tilde{X}$. In other words, it is convenient to introduce the chain $(X_p)_{0 \leq p \leq n} = (\tilde{X}_{kp})_{0 \leq p \leq n}$ where $k = \tilde{n}/n$ and $n$ is in the range 10–100. The underlying process can be a time-continuous Markov process $(\tilde{X}_t)_{t \in [0,T]}$ as well. In such a situation it is convenient to consider the chain $(X_p)_{0 \leq p \leq n} = (\tilde{X}_{pT/n})_{0 \leq p \leq n}$.

Beside the modeling of a new particle methodology, our main contribution is to provide a detailed asymptotic study of particle approximation models. Following the analysis of local sampling errors introduced in Chapter 9 in the research monograph [4], we first obtain an asymptotic expansion of the bias introduced by the interaction mechanism. We also design an original fluctuation analysis of polynomial functions of particle random fields, to derive new central limit theorems for weighted genealogical tree-based occupation measures. The magnitude of the asymptotic variances and comparisons with traditional Monte Carlo strategies are discussed in the context of Gaussian models.

Briefly, the paper is organized as follows. Section 2 contains all the theoretical results formulated in an abstract framework. We give a summary of the method and present a user-friendly implementation in Section 3. We consider a toy model (a Gaussian random walk) in Section 4 to illustrate the theoretical predictions on an example where all relevant quantities can be explicitly computed. Finally, in Section 5, we apply the method to a physical situation emerging in telecommunication.

## 2. Simulations of rare events by interacting particle systems.

2.1. *Introduction.* In this section we design an original IPS approach for analyzing Markov chains evolving in a rare event regime.

In Section 2.2 we use a natural large deviation perspective to exhibit natural changes of reference measures under which the underlying process is more likely to enter in a given rare level set. This technique is more or less well known. It often offers a powerful and elegant strategy for analyzing rare deviation probabilities. Loosely speaking, the twisted distributions associated to the deviated process represent the evolution of the original process in the rare event regime. In MC Markov chain literature, this changes-of-measure strategy is also called the importance sampling (IS) technique.

In Section 2.3 we present a Feynman–Kac formulation of twisted reference path distributions. We examine a pair of Gaussian models for which these changes of measures have a nice explicit formulation. In this context, we initiate a comparison of the fluctuation-error variances of the "pure" MC and the IS techniques. In general, the twisted distribution suggested by the physical model is rather complex, and its numerical analysis often requires extensive calculations. The practitioners often need to resort to another "suboptimal" reference strategy, based on a more refined analysis of the physical problem at hand. The main object of this section is to complement this IS methodology, by presenting a genetic type particle interpretation of a general and abstract class of twisted path models. Instead of hand crafting or simplified simulation models, this new particle methodology provides a powerful and very flexible way to produce samples according to any complex twisted measures dictated by the physical properties of the model at hand. But, from the strict practical point of view, if there exists already a good specialized IS method for a specific rare event problem, then our IPS methodology may not be the best tool for that application.

In Section 2.4 we introduce the reader to a new developing genealogical tree interpretation of Feynman–Kac path measures. For a more thorough study on this theme we refer to the monograph [4] and references therein. We connect this IPS methodology with rare event analysis. Intuitively speaking, the ancestral lines associated to these genetic evolution models represent the physical ways that the process uses to reach the desired rare level set.

In the final Section 2.5 we analyze the fluctuations of rare event particle simulation models. We discuss the performance of these interpretations on a class of warm-up Gaussian models. We compare the asymptotic error-variances of genealogical particle models and the more traditional noninteracting IS schemes. For Gaussian models, we show that the exponential fluctuation orders between these two particle simulation strategies are equivalent.

2.2. *A large deviation perspective.* Let $X_n$ be a Markov chain taking values at each time $n$ in some measurable state space $(E_n, \mathcal{E}_n)$ that may depend on the time parameter $n$. Suppose that we want to estimate the probability $P_n(a)$ that $X_n$ enters, at a given fixed date $n$, into the $a$-level set $V_n^{-1}([a, \infty))$ of a given energy-like function $V_n$ on $E_n$, for some $a \in \mathbb{R}$:

(2.1) $$P_n(a) = \mathbb{P}(V_n(X_n) \geq a).$$

To avoid some unnecessary technical difficulties, we further assume that $P_n(a) > 0$, and the pair $(X_n, V_n)$ satisfies Cramér's condition $\mathbb{E}(e^{\lambda V_n(X_n)}) < \infty$ for all $\lambda \in \mathbb{R}$. This condition ensures the exponential decay of the probabilities $\mathbb{P}(V_n(X_n) \geq a) \downarrow 0$, as $a \uparrow \infty$. To see this claim, we simply use the exponential version of Chebyshev's inequality to check that, for any $\lambda > 0$ we have

$$\mathbb{P}(V_n(X_n) \geq a) \leq e^{-\lambda(a - \lambda^{-1}\Lambda_n(\lambda))} \qquad \text{with } \Lambda_n(\lambda) \stackrel{\text{def.}}{=} \log \mathbb{E}(e^{\lambda V_n(X_n)}).$$

As an aside, it is also routine to prove that the maximum of $(\lambda a - \Lambda_n(\lambda))$ with respect to the parameter $\lambda > 0$ is attained at the value $\lambda_n(a)$ determined by the equation $a = \mathbb{E}(V_n(X_n)e^{\lambda V_n(X_n)}))/\mathbb{E}(e^{\lambda V_n(X_n)}))$. The resulting inequality

$$\mathbb{P}(V_n(X_n) \geq a) \leq e^{-\Lambda_n^\star(a)} \qquad \text{with } \Lambda_n^\star(a) = \sup_{\lambda > 0}(\lambda a - \Lambda_n(\lambda))$$

is known as large deviation inequality. When the Laplace transforms $\Lambda_n$ are explicitly known, this variational analysis often provides sharp tail estimates. We illustrate this observation on an elementary Gaussian model. This warm-up example will be used in several places in the further development of this article. In the subsequent analysis, it is briefly used primarily to carry out some variance calculations for natural IS strategies. As we already mentioned in the Introduction, and in this Gaussian context, we shall derive in Section 2.7 sharp estimates of mean error variances associated to a pair of IPS approximation models.

Suppose that $X_n$ is given by the recursive equation

(2.2) $$X_p = X_{p-1} + W_p$$

where $X_0 = 0$ and $(W_p)_{p \in \mathbb{N}^*}$ represents a sequence of independent and identically distributed (i.i.d.) Gaussian random variables, with $(\mathbb{E}(W_1), \mathbb{E}(W_1^2)) = (0,1)$. If we take $V_n(x) = x$, then we find that $\Lambda_n(\lambda) = \lambda^2 n/2$, $\lambda_n(a) = a/n$ and $\Lambda_n^\star(a) = a^2/(2n)$, from which we recover the well-known sharp exponential tails $\mathbb{P}(X_n \geq a) \leq e^{-a^2/(2n)}$.

In more general situations, the analytical expression of $\Lambda_n^\star(a)$ is out of reach, and we need to resort to judicious numerical strategies. The first rather crude MC method is to consider the estimate

$$P_n^N(a) = \frac{1}{N} \sum_{i=1}^N \mathbf{1}_{V_n(X_n^i) \geq a}$$

based on $N$ independent copies $(X_n^i)_{1 \leq i \leq N}$ of $X_n$. If is not difficult to check that the resulting error-variance is given by

$$\sigma_n^2(a) = N\mathbb{E}[(P_n^N(a) - P_n(a))^2] = P_n(a)(1 - P_n(a)).$$

In practice, $P_n^N(a)$ is a very poor estimate mainly because the whole sample set is very unlikely to reach the rare level.

A more judicious choice of MC exploration model is dictated by the large deviation analysis presented above. To be more precise, let us suppose that $a > \lambda^{-1}\Lambda_n(\lambda)$, with $\lambda > 0$. To simplify the presentation, we also assume that the initial value $X_0 = x_0$ is fixed, and we set $V_0(x_0) = 0$. Let $\mathbb{P}_n^\lambda$ be the new reference measure on the path space $F_n \stackrel{\mathrm{def.}}{=} (E_0 \times \cdots \times E_n)$ defined by the formula

$$(2.3) \qquad d\mathbb{P}_n^{(\lambda)} = \frac{1}{\mathbb{E}(e^{\lambda V_n(X_n)})} e^{\lambda V_n(X_n)}\, d\mathbb{P}_n,$$

where $\mathbb{P}_n$ is the distribution of the original and canonical path $(X_p)_{0 \leq p \leq n}$. By construction, we have that

$$\begin{aligned}\mathbb{P}(V_n(X_n) \geq a) &= \mathbb{E}_n^{(\lambda)}[\mathbf{1}_{V_n(X_n) \geq a}\, d\mathbb{P}_n/d\mathbb{P}_n^{(\lambda)}] \\ &= \mathbb{E}_n^{(\lambda)}[\mathbf{1}_{V_n(X_n) \geq a}\, e^{-\lambda V_n(X_n)}]\mathbb{E}[e^{\lambda V_n(X_n)}] \\ &\leq e^{-\lambda(a-\lambda^{-1}\Lambda_n(\lambda))}\mathbb{P}_n^{(\lambda)}(V_n(X_n) \geq a),\end{aligned}$$

where $\mathbb{E}_n^{(\lambda)}$ represents the expectation operator with respect to the distribution $\mathbb{P}_n^{(\lambda)}$. By definition, the measure $\mathbb{P}_n^{(\lambda)}$ tends to favor random evolutions with high potential values $V_n(X_n)$. As a consequence, the random paths under $\mathbb{P}_n^{(\lambda)}$ are much more likely to enter into the rare level set. For instance, in the Gaussian example described earlier, we have that

$$(2.4) \qquad d\mathbb{P}_n^{(\lambda)}/d\mathbb{P}_n = \prod_{p=1}^n e^{\lambda(X_p - X_{p-1}) - \lambda^2/2}.$$

In other words, under $\mathbb{P}_n^{(\lambda)}$ the chain takes the form $X_p = X_{p-1} + \lambda + W_p$, and we have $\mathbb{P}_n^{(\lambda)}(X_n \geq a) = \mathbb{P}_n(X_n \geq a - \lambda n)$ ($= 1/2$ as soon as $a = \lambda n$).

These observations suggest to replace $P^N(a)$ by the weighted MC model

$$P_n^{N,\lambda}(a) = \frac{1}{N}\sum_{i=1}^N \frac{d\mathbb{P}_n}{d\mathbb{P}_n^{(\lambda)}}(X_0^{\lambda,i},\ldots,X_n^{\lambda,i})\mathbf{1}_{V_n(X_n^{\lambda,i}) \geq a}$$

associated to $N$ independent copies $(X_n^{\lambda,i})_{1 \leq i \leq N}$ of the chain under $\mathbb{P}_n^{(\lambda)}$. Observe that the corresponding error-variance is given by

$$(2.5) \qquad \begin{aligned}\sigma_n^{(\lambda)}(a)^2 &= N\mathbb{E}[(P_n^{N,\lambda}(a) - P_n(a))^2] \\ &= \mathbb{E}[\mathbf{1}_{V_n(X_n) \geq a}\, e^{-\lambda V_n(X_n)}]\mathbb{E}[e^{\lambda V_n(X_n)}] - P_n^2(a) \\ &\leq e^{-\lambda(a-\lambda^{-1}\Lambda_n(\lambda))} P_n(a) - P_n^2(a).\end{aligned}$$

For judicious choices of $\lambda$, one expects the exponential large deviation term to be proportional to the desired tail probabilities $P_n(a)$. In this case, we have $\sigma_n^{(\lambda)}(a)^2 \leq KP_n^2(a)$ for some constant $K$. Returning to the Gaussian situation, and using Mill's inequalities

$$\frac{1}{t+1/t} \leq \mathbb{P}(\mathcal{N}(0,1) \geq t)\sqrt{2\pi}e^{t^2/2} \leq \frac{1}{t}$$

which are valid for any $t > 0$, and any reduced Gaussian random variable $\mathcal{N}(0,1)$ (see, e.g., (6) on page 237 in [25]), we find that

$$\sigma_n^{(\lambda)}(a)^2 \leq e^{-a^2/(2n)}P_n(a) - P_n^2(a) \leq P_n^2(a)[\sqrt{2\pi}(a/\sqrt{n} + \sqrt{n}/a) - 1]$$

for the value $\lambda = \lambda_n(a) = a/n$ which optimizes the large deviation inequality (2.6). For typical Gaussian type level indexes $a = a_0\sqrt{n}$, with large values of $a_0$, we find that $\lambda_n(a) = a_0/\sqrt{n}$ and

$$\sigma_n^{(\lambda)}(a)^2 \leq P_n^2(a)[\sqrt{2\pi}(a_0 + 1/a_0) - 1].$$

As an aside, although we shall be using most of the time upper bound estimates, Mill's inequalities ensure that most of the Gaussian exponential deviations are sharp.

The formulation (2.5) also suggests a dual interpretation of the variance. First, we note that

$$d\mathbb{P}_n/d\mathbb{P}_n^{(\lambda)} = \mathbb{E}[e^{\lambda V_n(X_n)}]\mathbb{E}[e^{-\lambda V_n(X_n)}]\, d\mathbb{P}_n^{(-\lambda)}/d\mathbb{P}_n$$

and therefore

$$\sigma_n^{(\lambda)}(a)^2 = \mathbb{P}_n^{(-\lambda)}(V_n(X_n) \geq a)\mathbb{E}[e^{\lambda V_n(X_n)}]\mathbb{E}[e^{-\lambda V_n(X_n)}] - P_n^2(a).$$

In contrast to $\mathbb{P}_n^{(\lambda)}$, the measure $\mathbb{P}_n^{(-\lambda)}$ now tends to favor low energy states $X_n$. As a consequence, we expect $\mathbb{P}_n^{(-\lambda)}(V_n(X_n) \geq a)$ to be much smaller than $P_n(a)$. For instance, in the Gaussian case we have

$$\mathbb{P}_n^{(-\lambda)}(X_n \geq a) = \mathbb{P}_n(X_n \geq a + \lambda n) \leq e^{-(a+\lambda n)^2/(2n)}.$$

Since we have $\mathbb{E}[e^{\lambda X_n}] = \mathbb{E}[e^{-\lambda X_n}] = e^{\lambda^2 n/2}$, we can write

(2.6) $$\sigma_n^{(\lambda)}(a)^2 \leq e^{-a^2/n}e^{(a-\lambda n)^2/(2n)} - P_n^2(a)$$

which confirms that the optimal choice (giving rise to the minimal variance) for the parameter $\lambda$ is $\lambda = a/n$.

2.3. *Twisted Feynman–Kac path measures.* The choice of the "twisted" measures $\mathbb{P}_n^{(\lambda)}$ introduced in (2.3) is only of pure mathematical interest. Indeed, the IS estimates described below will still require both the sampling of random paths according to $\mathbb{P}_n^{(\lambda)}$ and the computation of the normalizing constants. As we mentioned in the Introduction, the key difficulty in applying IS strategies is to choose the so-called "twisted" reference measures. In the further development of Section 2.4, we shall present a natural genealogical tree-based simulation technique of twisted Feynman–Kac path distribution of the following form:

$$(2.7) \qquad d\mathbb{Q}_n = \frac{1}{\mathcal{Z}_n} \left\{ \prod_{p=1}^{n} G_p(X_0, \ldots, X_p) \right\} d\mathbb{P}_n.$$

In the above display, $\mathcal{Z}_n > 0$ stands for a normalizing constant, and $(G_p)_{1 \leq p \leq n}$ represents a given sequence of potential functions on the path spaces $(F_p)_{1 \leq p \leq n}$. Note that the twisted measures defined in (2.3) correspond to the (nonunique) choice of functions

$$(2.8) \qquad G_p(X_0, \ldots, X_p) = e^{\lambda(V_p(X_p) - V_{p-1}(X_{p-1}))}.$$

As an aside, we mention that the optimal choice of twisted measure with respect to the IS criterion is the one associated to the potential functions $G_n = \mathbf{1}_{V_n^{-1}([a, \infty))}$, and $G_p = 1$, for $p < n$. In this case, we have $\mathcal{Z}_n = \mathbb{P}(V_n(X_n) \geq a)$ and $\mathbb{Q}_n$ is the distribution of the random paths ending in the desired rare level. This optimal choice is clearly infeasible, but we note that the resulting variance is null.

The rare event probability admits the following elementary Feynman–Kac formulation.

$$\mathbb{P}(V_n(X_n) \geq a) = \mathbb{E}\left[ g_n^{(a)}(X_0, \ldots, X_n) \prod_{p=1}^{n} G_p(X_0, \ldots, X_p) \right] = \mathcal{Z}_n \mathbb{Q}_n(g_n^{(a)})$$

with the weighted function defined by

$$g_n^{(a)}(x_0, \ldots, x_n) = \mathbf{1}_{V_n(x_n) \geq a} \prod_{p=1}^{n} G_p^{-1}(x_0, \ldots, x_p)$$

for any path sequence such that $\prod_{p=1}^{n} G_p(x_0, \ldots, x_p) > 0$. Otherwise, $g_n^{(a)}$ is assumed to be null.

The discussion given above already shows the improvements one might expect in changing the reference exploration measure. The central idea behind this IS methodology is to choose a twisted probability that mimics the physical behavior of the process in the rare event regime. The potential functions $G_p$ represent the changes of probability mass, and in some sense the physical variations in the evolution of the process to the rare level set. For instance, for time-homogeneous models $V_p = V$, $0 \leq p \leq n$, the potential functions defined in (2.8) will tend to favor local transitions that increase a given $V$-energy function. The large deviation analysis developed in Section 2.2 combined with the Feynman–Kac formulation (2.3) gives some indications on the way to choose the twisted potential functions $(G_p)_{1 \leq p \leq n}$.

Intuitively, the attractive forces induced by a particular choice of potentials are compensated by increasing normalizing constants. More formally, the error-variance of the $\mathbb{Q}_n$-importance sampling scheme is given by the formula

$$(2.9) \qquad \sigma_n^{\mathbb{Q}}(a)^2 = \mathbb{Q}_n^-(V_n(X_n) \geq a)\mathcal{Z}_n\mathcal{Z}_n^- - P_n(a)^2,$$

where $\mathbb{Q}_n^-$ is the path Feynman–Kac measure given by

$$d\mathbb{Q}_n^- = \frac{1}{\mathcal{Z}_n^-}\left\{\prod_{p=1}^n G_p^{-1}(X_0,\ldots,X_p)\right\} d\mathbb{P}_n.$$

Arguing as before, and since $\mathbb{Q}_n^-$ tends to favor random paths with low $G_p$ energy, we expect $\mathbb{Q}_n^-(V_n(X_n) \geq a)$ to be much smaller than the rare event probability $\mathbb{P}(V_n(X_n) \geq a)$. On the other hand, by Jensen's inequality we expect the product of normalizing constants $\mathcal{Z}_n\mathcal{Z}_n^-$ ($\geq 1$) to be very large. These expectations fail in the "optimal" situation examined above ($G_n = \mathbf{1}_{V_n^{-1}([a,\infty))}$, and $G_p = 1$, for $p < n$). In this case, we simply note that $\mathbb{Q}_n = \mathbb{Q}_n^- = \text{Law}(X_0,\ldots,X_n|V_n(X_n) \geq a)$, and $\mathbb{Q}_n^-(V_n(X_n) \geq a) = 1$, and $\mathcal{Z}_n = \mathcal{Z}_n^- = P_n(a)$. To avoid some unnecessary technical discussions, we always implicitly assume that the rare event probabilities $P_n(a)$ are strictly positive, so that the normalizing constants $\mathcal{Z}_n = \mathcal{Z}_n^- = P_n(a) > 0$ are always well defined.

We end this section with a brief discussion on the competition between making a rare event more attractive and controlling the normalizing constants. We return to the Gaussian example examined in (2.2), and instead of (2.4), we consider the twisted measure

$$(2.10) \qquad d\mathbb{Q}_n = d\mathbb{P}_n^{(\lambda)} = \frac{1}{\mathcal{Z}_n^{(\lambda)}}\left\{\prod_{p=1}^n e^{\lambda X_p}\right\} d\mathbb{P}_n.$$

In this case, it is not difficult to check that for any $\lambda \in \mathbb{R}$ we have $\mathcal{Z}_n^{(\lambda)} = e^{(\lambda^2/2)\sum_{p=1}^n p^2}$. In addition, under $\mathbb{P}_n^{(\lambda)}$ the chain $X_n$ has the form

$$(2.11) \qquad X_p = X_{p-1} + \lambda(n-p+1) + W_p, \qquad 1 \leq p \leq n.$$

When $\lambda > 0$, the rare level set is now very attractive, but the normalizing constants can become very large $\mathcal{Z}_n^{(\lambda)} = \mathcal{Z}_n^{(-\lambda)}$ ($\geq e^{\lambda^2 n^3/12}$). Also notice that in this situation the first term in the right-hand side of (2.9) is given by

$$\mathbb{P}_n^{(-\lambda)}(V_n(X_n) \geq a)\mathcal{Z}_n^{(\lambda)}\mathcal{Z}_n^{(-\lambda)}$$
$$\leq e^{(-1/(2n))(a+\lambda\sum_{p=1}^n p)^2 + \lambda^2 \sum_{p=1}^n p^2}$$
$$\leq e^{-a^2/n} e^{(1/(2n))(a-\lambda\sum_{p=1}^n p)^2 + \lambda^2[\sum_{p=1}^n p^2 - (\sum_{p=1}^n p)^2/n]}.$$

Although we are using inequalities, we recall that these exponential estimates are sharp. Now, if we take $\lambda = 2a/[n(n+1)]$, then we find that

$$(2.12) \qquad \mathbb{P}_n^{(-\lambda)}(V_n(X_n) \geq a)\mathcal{Z}_n^{(\lambda)}\mathcal{Z}_n^{(-\lambda)} \leq e^{-(a^2/n)(2/3)(1+1/(n+1))}.$$

This shows that even if we adjust correctly the parameter $\lambda$, this IS estimate is less efficient than the one associated to the twisted distribution (2.4). The reader has probably noticed that the change of measure defined in (2.10) is more adapted to estimate the probability of the rare level sets $\{V_n(Y_n) \geq a\}$, with the historical chain $Y_n = (X_0, \ldots, X_n)$ and the energy function $V_n(Y_n) = \sum_{p=1}^n X_p$.

2.4. *A genealogical tree-based interpretation model.* The probabilistic interpretation of the twisted Feynman–Kac measures (2.7) presented in this section can be interpreted as a mean field path-particle approximation of the distribution flow $(\mathbb{Q}_n)_{n \geq 1}$. We also mention that the genetic type selection/mutation evolution of the former algorithm can also be seen as an acceptance/rejection particle simulation technique. In this connection, and as we already mentioned in the Introduction, we again emphasize that this IPS methodology is not useful if we already know a specialized and exact simulation technique of the desired twisted measure.

2.4.1. *Rare event Feynman–Kac type distributions.* To simplify the presentation, it is convenient to formulate these models in terms of the historical process

$$Y_n \stackrel{\text{def.}}{=} (X_0, \ldots, X_n) \in F_n \stackrel{\text{def.}}{=} (E_0 \times \cdots \times E_n).$$

We let $M_n(y_{n-1}, dy_n)$ be the Markov transitions associated to the chain $Y_n$. To simplify the presentation, we assume that the initial value $Y_0 = X_0 = x_0$ is fixed, and we also denote by $K_n(x_{n-1}, dx_n)$ the Markov transitions of $X_n$. We finally let $\mathcal{B}_b(E)$ be the space of all bounded measurable functions on some measurable space $(E, \mathcal{E})$, and we equip $\mathcal{B}_b(E)$ with the uniform norm.

We associate to the pair potentials/transitions $(G_n, M_n)$ the Feynman–Kac measure defined for any test function $f_n \in \mathcal{B}_b(F_n)$ by the formula

$$\gamma_n(f_n) = \mathbb{E}\left[f_n(Y_n) \prod_{1 \leq k < n} G_k(Y_k)\right].$$

We also introduce the corresponding normalized measure

$$\eta_n(f_n) = \gamma_n(f_n)/\gamma_n(1).$$

To simplify the presentation and avoid unnecessary technical discussions, we suppose that the potential functions are chosen such that

$$\sup_{(y_n, y_n') \in F_n^2} G_n(y_n)/G_n(y_n') < \infty.$$

This regularity condition ensures that the normalizing constants $\gamma_n(1)$ and the measure $\gamma_n$ are bounded and positive. This technical assumption clearly fails for unbounded or for indicator type potential functions. The Feynman–Kac and the particle approximation models developed in this section can be extended to more general situations using traditional cut-off techniques, by considering Kato-class type of potential functions (see, e.g., [19, 22, 26]), or by using different Feynman–Kac representations of the twisted measures (see, e.g., Section 2.5 in [4]).

In this section we provide a Feynman–Kac formulation of rare event probabilities. The fluctuation analysis of their genealogical tree interpretations will also be described in terms of the distribution flow $(\gamma_n^-, \eta_n^-)$, defined as $(\gamma_n, \eta_n)$ by replacing the potential functions $G_p$ by their inverse

$$G_p^- = 1/G_p.$$

The twisted measures $\mathbb{Q}_n$ presented in (2.7) and the desired rare event probabilities have the following Feynman–Kac representation:

$$\mathbb{Q}_n(f_n) = \eta_n(f_n G_n)/\eta_n(G_n) \quad \text{and} \quad \mathbb{P}(V_n(X_n) \geq a) = \gamma_n(T_n^{(a)}(1)).$$

In the above displayed formulae, $T_n^{(a)}(1)$ is the weighted indicator function defined for any path $y_n = (x_0, \ldots, x_n) \in F_n$ by

$$T_n^{(a)}(1)(y_n) = T_n^{(a)}(1)(x_0, \ldots, x_n) = \mathbf{1}_{V_n(x_n) \geq a} \prod_{1 \leq p < n} G_p^-(x_0, \ldots, x_p).$$

More generally, we have for any $\varphi_n \in \mathcal{B}_b(F_n)$

$$\mathbb{E}[\varphi_n(X_0, \ldots, X_n); V_n(X_n) \geq a] = \gamma_n(T_n^{(a)}(\varphi_n))$$

with the function $T_n^{(a)}(\varphi_n)$ given by

$$(2.13) \quad T_n^{(a)}(\varphi_n)(x_0, \ldots, x_n) = \varphi_n(x_0, \ldots, x_n) \mathbf{1}_{V_n(x_n) \geq a} \prod_{1 \leq p < n} G_p^-(x_0, \ldots, x_p).$$

To connect the rare event probabilities with the normalized twisted measures we use the fact that

$$\gamma_{n+1}(1) = \gamma_n(G_n) = \eta_n(G_n)\gamma_n(1) = \prod_{p=1}^n \eta_p(G_p).$$

This readily implies that for any $f_n \in \mathcal{B}_b(F_n)$

$$(2.14) \quad \gamma_n(f_n) = \eta_n(f_n) \prod_{1 \leq p < n} \eta_p(G_p).$$

This yields the formulae

$$\mathbb{P}(V_n(X_n) \geq a) = \eta_n(T_n^{(a)}(1)) \prod_{1 \leq p < n} \eta_p(G_p),$$

(2.15) $\quad \mathbb{E}(\varphi_n(X_0, \ldots, X_n); V_n(X_n) \geq a) = \eta_n(T_n^{(a)}(\varphi_n)) \prod_{1 \leq p < n} \eta_p(G_p),$

$$\mathbb{E}(\varphi_n(X_0, \ldots, X_n) | V_n(X_n) \geq a) = \eta_n(T_n^{(a)}(\varphi_n))/\eta_n(T_n^{(a)}(1)).$$

To take the final step, we use the Markov property to check that the twisted measures $(\eta_n)_{n \geq 1}$ satisfy the nonlinear recursive equation

$$\eta_n = \Phi_n(\eta_{n-1}) \stackrel{\text{def.}}{=} \int_{F_{n-1}} \eta_{n-1}(dy_{n-1}) G_{n-1}(y_{n-1}) M_n(y_{n-1}, \cdot)/\eta_{n-1}(G_{n-1})$$

starting from $\eta_1 = M_1(x_0, \cdot)$.

2.4.2. *Interacting path-particle interpretation.* The mean field particle model associated with a collection of transformations $\Phi_n$ is a Markov chain $\xi_n = (\xi_n^i)_{1 \leq i \leq N}$ taking values at each time $n \geq 1$ in the product spaces $F_n^N$. Loosely speaking, the algorithm will be conducted so that each path-particle

$$\xi_n^i = (\xi_{0,n}^i, \xi_{1,n}^i, \ldots, \xi_{n,n}^i) \in F_n = (E_0 \times \cdots \times E_n)$$

is almost sampled according to the twisted measure $\eta_n$.

The initial configuration $\xi_1 = (\xi_1^i)_{1 \leq i \leq N}$ consists of $N$ independent and identically distributed random variables with common distribution

$$\eta_1(d(y_0, y_1)) = M_1(x_0, d(y_0, y_1)) = \delta_{x_0}(dy_0) \ K_1(y_0, dy_1).$$

In other words, $\xi_1^i \stackrel{\text{def.}}{=} (\xi_{0,1}^i, \xi_{1,1}^i) = (x_0, \xi_{1,1}^i) \in F_1 = (E_0 \times E_1)$ can be interpreted as $N$ independent copies $x_0 \rightsquigarrow \xi_{1,1}^i$ of the initial elementary transition $X_0 = x_0 \rightsquigarrow X_1$. The elementary transitions $\xi_{n-1} \rightsquigarrow \xi_n$ from $F_{n-1}^N$ into $F_n^N$ are defined by

(2.16) $\quad \mathbb{P}(\xi_n \in d(y_n^1, \ldots, y_n^N) | \xi_{n-1}) = \prod_{j=1}^N \Phi_n(m(\xi_{n-1}))(dy_n^j),$

where $m(\xi_{n-1}) \stackrel{\text{def.}}{=} \frac{1}{N} \sum_{i=1}^N \delta_{\xi_{n-1}^i}$, and $d(y_n^1, \ldots, y_n^N)$ is an infinitesimal neighborhood of the point $(y_n^1, \ldots, y_n^N) \in F_n^N$. By the definition of $\Phi_n$ we find that (2.16) is the overlapping of simple selection/mutation genetic transitions

$$\xi_{n-1} \in F_{n-1}^N \stackrel{\text{selection}}{\longrightarrow} \hat{\xi}_{n-1} \in F_{n-1}^N \stackrel{\text{mutation}}{\longrightarrow} \xi_n \in F_n^N.$$

The selection stage consists of choosing randomly and independently $N$ path-particles

$$\hat{\xi}^i_{n-1} = (\hat{\xi}^i_{0,n-1}, \hat{\xi}^i_{1,n-1}, \ldots, \hat{\xi}^i_{n-1,n-1}) \in F_{n-1}$$

according to the Boltzmann–Gibbs particle measure

$$\sum_{j=1}^{N} \frac{G_{n-1}(\xi^j_{0,n-1}, \ldots, \xi^j_{n-1,n-1})}{\sum_{j'=1}^{N} G_{n-1}(\xi^{j'}_{0,n-1}, \ldots, \xi^{j'}_{n-1,n-1})} \, \delta_{(\xi^j_{0,n-1}, \ldots, \xi^j_{n-1,n-1})}.$$

During the mutation stage, each selected path-particle $\hat{\xi}^i_{n-1}$ is extended by an elementary $K_n$-transition. In other words, we set

$$\xi^i_n = ((\xi^i_{0,n}, \ldots, \xi^i_{n-1,n}), \xi^i_{n,n})$$
$$= ((\hat{\xi}^i_{0,n-1}, \ldots, \hat{\xi}^i_{n-1,n-1}), \xi^i_{n,n}) \in F_n = F_{n-1} \times E_n,$$

where $\xi^i_{n,n}$ is a random variable with distribution $K_n(\hat{\xi}^i_{n-1,n-1}, \cdot)$. The mutations are performed independently.

2.4.3. *Particle approximation measures.* It is of course out of the scope of this article to present a full asymptotic analysis of these genealogical particle models. We rather refer the interested reader to the recent monograph [4], and the references therein. For instance, it is well known that the occupation measures of the ancestral lines converge to the desired twisted measures. That is, we have with various precision estimates the weak convergence result

(2.17) $$\eta_n^N \stackrel{\text{def.}}{=} \frac{1}{N} \sum_{i=1}^{N} \delta_{(\xi^i_{0,n}, \xi^i_{1,n}, \ldots, \xi^i_{n,n})} \stackrel{N \to \infty}{\longrightarrow} \eta_n.$$

In addition, several propagation-of-chaos estimates ensure that the ancestral lines $(\xi^i_{0,n}, \xi^i_{1,n}, \ldots, \xi^i_{n,n})$ are asymptotically independent and identically distributed with common distribution $\eta_n$. The asymptotic analysis of regular models with unbounded potential functions can be treated using traditional cut-off techniques.

Mimicking (2.14), the *unbias* particle approximation measures $\gamma_n^N$ of the unnormalized model $\gamma_n$ are defined as

$$\gamma_n^N(f_n) = \eta_n^N(f_n) \prod_{1 \leq p < n} \eta_p^N(G_p).$$

By (2.15), we eventually get the particle approximation of the rare event probabilities $P_n(a)$. More precisely, if we let

(2.18) $$P_n^N(a) = \gamma_n^N(T_n^{(a)}(1)) = \eta_n^N(T_n^{(a)}(1)) \prod_{1 \leq p < n} \eta_p^N(G_p),$$

then we find that $P_n^N(a)$ is an unbiased estimator of $P_n(a)$ such that

$$P_n^N(a) \overset{N \to \infty}{\longrightarrow} P_n(a) \qquad \text{a.s.} \tag{2.19}$$

In addition, by (2.15), the conditional distribution of the process in the rare event regime can be estimated using the weighted particle measure

$$\begin{aligned}
P_n^N(a, \varphi_n) &\overset{\text{def.}}{=} \eta_n^N(T_n^{(a)}(\varphi_n))/\eta_n^N(T_n^{(a)}(1)) \\
&\overset{N \to \infty}{\longrightarrow} P_n(a, \varphi_n) \overset{\text{def.}}{=} \mathbb{E}[\varphi_n(X_0, \ldots, X_n) | V_n(X_n) \geq a].
\end{aligned} \tag{2.20}$$

When no particles have succeeded in reaching the desired level $V_n^{-1}([a, \infty))$ at time $n$, we have $\eta_n^N(T_n^{(a)}(1)) = 0$, and therefore $\eta_n^N(T_n^{(a)}(\varphi_n)) = 0$ for any $\varphi_n \in \mathcal{B}_b(F_n)$. In this case, we take the convention $P_n^N(a, \varphi_n) = 0$. Also notice that $\eta_n^N(T_n^{(a)}(1)) > 0$ if and only if we have $P_n^N(a) > 0$. When $P_n(a) > 0$, we have the exponential decay of the probabilities $\mathbb{P}(P_n^N(a) = 0) \to 0$ as $N$ tends to infinity.

The above asymptotic, and reassuring, estimates are almost sure convergence results. Their complete proofs, together with the analysis of extinction probabilities, rely on a precise propagation-of-chaos type analysis. They can be found in Section 7.4, pages 239–241, and Theorem 7.4.1, page 232 in [4]. In Section 2.5 we provide a natural and simple proof of the consistency of the particle measures $(\gamma_n^N, \eta_n^N)$ using an original fluctuation analysis.

In our context, these almost sure convergence results show that the genealogical tree-based approximation schemes of rare event probabilities are consistent. Unfortunately, the rather crude estimates say little, as much as more naive numerical methods do converge as well. Therefore, we need to work harder to analyze the precise asymptotic bias and the fluctuations of the occupation measures of the complete and weighted genealogical tree. These questions, as well as comparisons of the asymptotic variances, are addressed in the next three sections.

We can already mention that the consistency results discussed above will be pivotal in the more refined analysis of particle random fields. They will be used in the further development of Section 2.5, in conjunction with a semigroup technique, to derive central limit theorems for particle random fields.

2.5. *Fluctuations analysis.* The fluctuations of genetic type particle models have been initiated in [5]. Under appropriate regularity conditions on the mutation transitions, this study provides a central limit theorem for the path-particle model $(\xi_0^i, \ldots, \xi_n^i)_{1 \leq i \leq N}$. Several extensions, including Donsker's type theorems, Berry–Esseen inequalities and applications to nonlinear filtering problems can be found in [6, 7, 8, 9]. In this section we design a simplified analysis essentially based on the fluctuations of random fields associated to the local sampling errors. In this subsection we provide a brief discussion on the fluctuations analysis of the weighted particle measures introduced in Section 2.4. We underline several interpretations of the central limit variances in terms of twisted Feynman–Kac measures. In the final part of this section we illustrate these general and theoretical fluctuations analyses in the warm-up Gaussian situation discussed in (2.2), (2.4) and (2.10). In this context, we derive an explicit description of the error-variances, and we compare the performance of the IPS methodology with the noninteracting IS technique.

The fluctuations of the mean field particle models described in Section 2.4 are essentially based on the asymptotic analysis of the local sampling errors associated with the particle approximation sampling steps. These local errors are defined in terms of the random fields $\mathcal{W}_n^N$, given for any $f_n \in \mathcal{B}_b(F_n)$ by the formula

$$\mathcal{W}_n^N(f_n) = \sqrt{N}\ [\eta_n^N - \Phi_n(\eta_{n-1}^N)](f_n).$$

The next central limit theorem is pivotal ([4], Theorem 9.3.1, page 295). For any fixed time horizon $n \geq 1$, the sequence $(\mathcal{W}_p^N)_{1 \leq p \leq n}$ converges in law, as $N$ tends to infinity, to a sequence of $n$ independent, Gaussian and centered random fields $(\mathcal{W}_p)_{1 \leq p \leq n}$; with, for any $f_p, g_p \in \mathcal{B}_b(F_p)$, and $1 \leq p \leq n$,

$$\mathbb{E}[\mathcal{W}_p(f_p)\mathcal{W}_p(g_p)] = \eta_p([f_p - \eta_p(f_p)][g_p - \eta_p(g_p)]).$$

Let $Q_{p,n}$, $1 \leq p \leq n$, be the Feynman–Kac semigroup associated to the distribution flow $(\gamma_p)_{1 \leq p \leq n}$. For $p = n$ it is defined by $Q_{n,n} = Id$, and for $p < n$ it has the following functional representation:

$$Q_{p,n}(f_n)(y_p) = \mathbb{E}\left[f_n(Y_n) \prod_{p \leq k < n} G_k(Y_k) | Y_p = y_p\right]$$

for any test function $f_n \in \mathcal{B}_b(F_n)$, and any path sequence $y_p = (x_0, \ldots, x_p) \in F_p$. The semigroup $Q_{p,n}$ satisfies

(2.21) $$\forall 1 \leq p \leq n \qquad \gamma_n = \gamma_p Q_{p,n}.$$

To check this assertion, we note that

$$\gamma_n(f_n) = \mathbb{E}\left[f_n(Y_n) \prod_{1 \leq k < n} G_k(Y_k)\right]$$

$$= \mathbb{E}\left[\left[\prod_{1 \leq k < p} G_k(Y_k)\right] \times \mathbb{E}\left(f_n(Y_n) \prod_{p \leq k < n} G_k(Y_k) | (Y_0, \ldots, Y_p)\right)\right]$$

for any $f_n \in \mathcal{B}_b(E_n)$. Using the Markov property we conclude that

$$\gamma_n(f_n) = \mathbb{E}\left[\left[\prod_{1 \leq k < p} G_k(Y_k)\right] \times \mathbb{E}\left(f_n(Y_n) \prod_{p \leq k < n} G_k(Y_k) | Y_p\right)\right]$$

$$= \mathbb{E}\left[\left[\prod_{1 \leq k < p} G_k(Y_k)\right] Q_{p,n}(f_n)(Y_p)\right] = \gamma_p Q_{p,n}(f_n)$$

which establishes (2.21). To explain what we have in mind when making these definitions, we now consider the elementary telescopic decomposition

$$\gamma_n^N - \gamma_n = \sum_{p=1}^n [\gamma_p^N Q_{p,n} - \gamma_{p-1}^N Q_{p-1,n}].$$

For $p = 1$, we recall that $\eta_0^N = \delta_{x_0}$ and $\gamma_1 = \eta_1 = M_1(x_0, \cdot)$, from which we find that $\eta_0^N Q_{0,n} = \gamma_1 Q_{1,n} = \gamma_n$. Using the fact that

$$\gamma_{p-1}^N Q_{p-1,p} = \gamma_{p-1}^N(G_{p-1}) \times \Phi_{p-1}(\eta_{p-1}^N) \quad \text{and} \quad \gamma_{p-1}^N(G_{p-1}) = \gamma_p^N(1)$$

the above decomposition implies that

$$(2.22) \quad \mathcal{W}_n^{\gamma,N}(f_n) \stackrel{\text{def.}}{=} \sqrt{N}[\gamma_n^N - \gamma_n](f_n) = \sum_{p=1}^n \gamma_p^N(1) \mathcal{W}_p^N(Q_{p,n} f_n).$$

LEMMA 2.1. $\gamma_n^N$ is an unbiased estimate of $\gamma_n$, in the sense that for any $p \geq 1$ and $f_n \in \mathcal{B}_b(F_n)$, with $\|f_n\| \leq 1$, we have

$$\mathbb{E}(\gamma_n^N(f_n)) = \gamma_n(f_n) \quad \text{and} \quad \sup_{N \geq 1} \sqrt{N} \mathbb{E}[|\gamma_n^N(f_n) - \gamma_n(f_n)|^p]^{1/p} \leq c_p(n)$$

for some constant $c_p(n) < \infty$ whose value does not depend on the function $f_n$.

PROOF. We first notice that $(\gamma_n^N(1))_{n \geq 1}$ is a predictable sequence, in the sense that

$$\mathbb{E}(\gamma_n^N(1) | \xi_0, \ldots, \xi_{n-1}) = \mathbb{E}\left(\prod_{p=1}^{n-1} \eta_p^N(G_p) | \xi_0, \ldots, \xi_{n-1}\right) = \prod_{p=1}^{n-1} \eta_p^N(G_p) = \gamma_n^N(1).$$

On the other hand, by definition of the particle scheme, for any $1 \leq p \leq n$, we also have that

$$\mathbb{E}(\mathcal{W}_p^N(Q_{p,n} f_n) | \xi_0, \ldots, \xi_{p-1})$$
$$= \sqrt{N} \mathbb{E}([\eta_p^N - \Phi_p(\eta_{p-1}^N)](Q_{p,n} f_n) | \xi_0, \ldots, \xi_{p-1}) = 0.$$

Combining these two observations, we find that

$$\mathbb{E}(\gamma_p^N(1)\mathcal{W}_p^N(Q_{p,n}f_n)|\xi_0,\ldots,\xi_{p-1}) = 0.$$

This yields that $\gamma_n^N$ is unbiased. In the same way, using the fact that the potential functions are bounded, we have for any $p \geq 1$ and $f_n \in \mathcal{B}_b(F_n)$, with $\|f_n\| \leq 1$,

$$\mathbb{E}[|[\gamma_n^N - \gamma_n](f_n)|^p]^{1/p} \leq \sum_{k=1}^n a_1(k)\mathbb{E}[|(\eta_p^N - \Phi_p(\eta_{p-1}^N))(Q_{p,n}f_n)|^p]^{1/p}$$

for some constant $a_1(k) < \infty$ which only depends on the time parameter. We recall that $\eta_n^N$ is the empirical measure associated with a collection of $N$ conditionally independent particles with common law $\Phi_p(\eta_{p-1}^N)$. The end of the proof is now a consequence of Burkholder's inequality. $\square$

Lemma 2.1 shows that the random sequence $(\gamma_p^N(1))_{1 \leq p \leq n}$ converges in probability, as $N$ tends to infinity, to the deterministic sequence $(\gamma_p(1))_{1 \leq p \leq n}$. An application of Slutsky's lemma now implies that the random fields $\mathcal{W}_n^{\gamma,N}$ converge in law, as $N$ tends to infinity, to the Gaussian random fields $\mathcal{W}_n^{\gamma}$ defined for any $f_n \in \mathcal{B}_b(F_n)$ by

$$(2.23) \qquad \mathcal{W}_n^{\gamma}(f_n) = \sum_{p=1}^n \gamma_p(1)\mathcal{W}_p(Q_{p,n}f_n).$$

In much the same way, the sequence of random fields

$$(2.24) \qquad \begin{aligned} \mathcal{W}_n^{\eta,N}(f_n) &\stackrel{\text{def.}}{=} \sqrt{N}[\eta_n^N - \eta_n](f_n) \\ &= \frac{\gamma_n(1)}{\gamma_n^N(1)} \times \mathcal{W}_n^{\gamma,N}\left(\frac{1}{\gamma_n(1)}(f_n - \eta_n(f_n))\right) \end{aligned}$$

converges in law, as $N$ tends to infinity, to the Gaussian random fields $\mathcal{W}_n^{\eta}$ defined for any $f_n \in \mathcal{B}_b(F_n)$ by

$$(2.25) \qquad \begin{aligned} \mathcal{W}_n^{\eta}(f_n) &= \mathcal{W}_n^{\gamma}\left(\frac{1}{\gamma_n(1)}(f_n - \eta_n(f_n))\right) \\ &= \sum_{p=1}^n \mathcal{W}_p\left(\frac{Q_{p,n}}{\eta_p Q_{p,n}(1)}(f_n - \eta_n(f_n))\right). \end{aligned}$$

The key decomposition (2.24) also appears to be useful to obtain $\mathbb{L}_p$-mean error bounds. More precisely, recalling that $\gamma_n(1)/\gamma_n^N(1)$ is a uniformly bounded sequence w.r.t. the population parameter $N \geq 1$, and using Lemma 2.1, we prove the following result.

LEMMA 2.2. *For any $p \geq 1$ and $f_n \in \mathcal{B}_b(F_n)$, with $\|f_n\| \leq 1$, we have*

$$\sup_{N \geq 1} \sqrt{N} \mathbb{E}[|\eta_n^N(f_n) - \eta_n(f_n)|^p]^{1/p} \leq c_p(n)$$

*for some constant $c_p(n) < \infty$ whose value does not depend on the function $f_n$.*

A consequence of the above fluctuations is a central limit theorem for the estimators $P_n^N(a)$ and $P_n^N(a, \phi_n)$ introduced in (2.18) and (2.20).

THEOREM 2.3. *The estimator $P_n^N(a)$ given by (2.18) is unbiased, and it satisfies the central limit theorem*

$$(2.26) \qquad \sqrt{N}[P_n^N(a) - P_n(a)] \overset{N \to \infty}{\longrightarrow} \mathcal{N}(0, \sigma_n^\gamma(a)^2)$$

*with the asymptotic variance*

$$(2.27) \qquad \sigma_n^\gamma(a)^2 = \sum_{p=1}^n [\gamma_p(1)\gamma_p^-(1)\eta_p^-(P_{p,n}(a)^2) - P_n(a)^2]$$

*and the collection of functions $P_{p,n}(a)$ defined by*

$$(2.28) \qquad x_p \in E_p \mapsto P_{p,n}(a)(x_p) = \mathbb{E}[\mathbf{1}_{V_n(X_n) \geq a} | X_p = x_p] \in [0, 1].$$

*In addition, for any $\varphi_n \in \mathcal{B}_b(F_n)$, the estimator $P_n^N(a, \varphi_n)$ given by (2.20) satisfies the central limit theorem*

$$(2.29) \qquad \sqrt{N}[P_n^N(a, \varphi_n) - P_n(a, \varphi_n)] \overset{N \to \infty}{\longrightarrow} \mathcal{N}(0, \sigma_n(a, \varphi_n)^2)$$

*with the asymptotic variance*

$$(2.30) \qquad \sigma_n(a, \varphi_n)^2 = P_n(a)^{-2} \sum_{p=1}^n \gamma_p(1)\gamma_p^-(1)\eta_p^-(P_{p,n}(a, \varphi_n)^2)$$

*and the collection of functions $P_{p,n}(a, \varphi_n) \in \mathcal{B}_b(F_p)$ defined by*

$$(2.31) \quad \begin{aligned} P_{p,n}(a, \varphi_n)(x_0, \ldots, x_p) &= \mathbb{E}[(\varphi_n(X_0, \ldots, X_n) - P_n(a, \varphi_n))\mathbf{1}_{V_n(X_n) \geq a} | \\ &\quad (X_0, \ldots, X_p) = (x_0, \ldots, x_p)]. \end{aligned}$$

PROOF. We first notice that

$$(2.32) \qquad \sqrt{N} \, [P_n^N(a) - P_n(a)] = \mathcal{W}_n^{\gamma, N}(T_n^{(a)}(1))$$

with the weighted function $T_n^{(a)}(1)$ introduced in (2.13). If we take $f_n = T_n^{(a)}(1)$ in (2.22) and (2.23), then we find that $\mathcal{W}_n^{\gamma, N}(T_n^{(a)}(1))$ converges in law, as $N$ tends to infinity, to a centered Gaussian random variable $\mathcal{W}_n^\gamma(T_n^{(a)}(1))$ with the variance

$$\sigma_n^\gamma(a)^2 \stackrel{\text{def.}}{=} \mathbb{E}(\mathcal{W}_n^\gamma(T_n^{(a)}(1))^2)$$
$$= \sum_{p=1}^n \gamma_p(1)^2 \eta_p([Q_{p,n}(T_n^{(a)}(1)) - \eta_p Q_{p,n}(T_n^{(a)}(1))]^2).$$

To have a more explicit description of $\sigma_n^\gamma(a)$ we notice that

$$Q_{p,n}(T_n^{(a)}(1))(x_0,\ldots,x_p)$$
$$= \left\{\prod_{1 \leq k < p} G_k(x_0,\ldots,x_k)^{-1}\right\} \mathbb{P}(V_n(X_n) \geq a | X_p = x_p).$$

By definition of $\eta_p$, we also find that

$$\eta_p(Q_{p,n}(T_n^{(a)}(1))) = \mathbb{P}(V_n(X_n) \geq a)/\gamma_p(1).$$

From these observations, we conclude that

(2.33)
$$\sigma_n^\gamma(a)^2 = \sum_{p=1}^n \left\{\gamma_p(1)\mathbb{E}\left[\prod_{1 \leq k < p} G_k^-(X_0,\ldots,X_k)\mathbb{E}(\mathbf{1}_{V_n(X_n) \geq a}|X_p)^2\right]\right.$$
$$\left. - P_n(a)^2\right\}.$$

This variance can be rewritten in terms of the distribution flow $(\eta_p^-)_{1 \leq p \leq n}$, since we have

$$\gamma_p^-(f_p) = \mathbb{E}\left[\left[\prod_{1 \leq k < p} G^-(X_0,\ldots,X_k)\right] \times f_p(X_p)\right] = \eta_p^-(f_p) \times \gamma_p^-(1)$$

for any $f_p \in \mathcal{B}_b(E_p)$, and any $1 \leq p \leq n$. Substituting into (2.33) yields (2.27).

Our next objective is to analyze the fluctuations of the particle conditional distributions of the process in the rare event regime defined in (2.20):

$$\sqrt{N}[P_n^N(a,\varphi_n) - P_n(a,\varphi_n)] = \frac{\eta_n T_n^{(a)}(1)}{\eta_n^N T_n^{(a)}(1)} \times \mathcal{W}_n^{\eta,N}\left(\frac{T_n^{(a)}}{\eta_n T_n^{(a)}(1)}(\varphi_n - P_n(a,\varphi_n))\right).$$

Using the same arguments as above, one proves that the sequence of random variables $\frac{\eta_n T_n^{(a)}(1)}{\eta_n^N T_n^{(a)}(1)}$ converges in probability, as $N \to \infty$, to 1. Therefore, letting

$$f_n = \frac{T_n^{(a)}}{\eta_n T_n^{(a)}(1)}(\varphi_n - P_n(a,\varphi_n))$$

in (2.24) and (2.25), and applying again Slutsky's lemma, we have the weak convergence

$$\sqrt{N}[P_n^N(a,\varphi_n) - P_n(a,\varphi_n)]\mathbf{1}_{P_n^N(a)>0}$$

(2.34)
$$\stackrel{N\to\infty}{\Longrightarrow} \mathcal{W}_n^\eta\left(\frac{T_n^{(a)}}{\eta_n T_n^{(a)}(1)}(\varphi_n - P_n(a,\varphi_n))\right).$$

The limit is a centered Gaussian random variable with the variance

$$\sigma_n(a,\varphi_n)^2 \stackrel{\text{def.}}{=} \mathbb{E}\left[\mathcal{W}_n^\eta\left(\frac{T_n^{(a)}}{\eta_n T_n^{(a)}(1)}(\varphi_n - P_n(a,\varphi_n))\right)^2\right].$$

Taking into account the definition of $\mathcal{W}_n^\eta$ and the identities $\eta_n T_n^{(a)}(1) = P_n(a)/\gamma_n(1)$ and $\eta_n T_n^{(a)}(\varphi_n - P_n(a,\varphi_n)) = 0$, we obtain

$$(2.35)\quad \sigma_n(a,\varphi_n)^2 = P_n(a)^{-2} \sum_{p=1}^n \gamma_p(1)\gamma_p([Q_{p,n}(T_n^{(a)}(\varphi_n - P_n(a,\varphi_n)))]^2).$$

To derive (2.30) from (2.35), we notice that

$$Q_{p,n}(T_n^{(a)}(\varphi_n - P_n(a,\varphi_n)))(x_0,\ldots,x_p)$$
$$= \mathbb{E}\left[\left(\prod_{p\leq k<n} G(X_0,\ldots,X_k)\right)\mathbf{1}_{V_n(X_n)\geq a}\left(\prod_{1\leq k<n} G^-(X_0,\ldots,X_k)\right)\right.$$
$$\left. \times (\varphi_n(X_0,\ldots,X_n) - P_n(a,\varphi_n))|(X_0,\ldots,X_p) = (x_0,\ldots,x_p)\right]$$
$$= \left(\prod_{1\leq k<p} G^-(x_0,\ldots,x_k)\right)$$
$$\times \mathbb{E}[(\varphi_n(X_0,\ldots,X_n) - P_n(a,\varphi_n))$$
$$\times \mathbf{1}_{V_n(X_n)\geq a}|(X_0,\ldots,X_p) = (x_0,\ldots,x_p)]$$
$$= \left(\prod_{1\leq k<p} G^-(x_0,\ldots,x_k)\right) \times P_{p,n}(a,\varphi_n)(x_0,\ldots,x_p)$$

from which we find that

$$\gamma_p([Q_{p,n}(T_n^{(a)}(\varphi_n - P_n(a,\varphi_n)))]^2)$$
$$= \mathbb{E}\left[\left(\prod_{1\leq k<p} G^-(X_0,\ldots,X_k)\right) \times P_{p,n}(a,\varphi_n)(X_0,\ldots,X_p)^2\right]$$
$$= \gamma_p^-(P_{p,n}(a,\varphi_n)^2) = \gamma_p^-(1) \times \eta_p^-(P_{p,n}(a,\varphi_n)^2). \qquad \square$$

Arguing as in the end of Section 2.2, we note that the measures $\eta_p^-$ tend to favor random paths with low $(G_k)_{1\leq k<p}$ potential values. Recalling that these potentials are chosen so as to represent the process evolution in the rare level set, we expect the quantities $\eta_p^-(P_{p,n}(a)^2)$ to be much smaller than $P_n(a)$. In the reverse angle, by Jensen's inequality we expect the normalizing constants products $\gamma_p(1)\gamma_p^-(1)$ to be rather large. We shall make precise these intuitive comments in the next section, with explicit calculations for the pair Gaussian models introduced in (2.4) and (2.10). We end the section by noting that

$$\sigma_n(a,\varphi_n)^2 \leq P_n(a)^{-2}\sum_{p=1}^n \gamma_p(1)\gamma_p^-(1)\eta_p^-(P_{p,n}(a)^2)$$

for any test function $\varphi_n$, with $\sup_{(y_n,y_n')\in F_p^2}|\varphi_n(y_n)-\varphi_n(y_n')|\leq 1$.

2.6. *On the weak negligible bias of genealogical models.* In this subsection we complete the fluctuation analysis developed in Section 2.5 with the study of the bias of the genealogical tree occupation measures $\eta_n^N$, and the corresponding weighted measures $P_n^N(a,\varphi_n)$ defined by (2.20). The forthcoming analysis also provides sharp estimates, and a precise asymptotic description of the law of a given particle ancestral line. In this sense, this study also completes the propagation-of-chaos analysis developed in [4]. The next technical lemma is pivotal in our way to analyze the bias of the path-particle models.

LEMMA 2.4. *For any $n,d\geq 1$, any collection of functions $(f_n^i)_{1\leq i\leq d}\in \mathcal{B}_b(F_n)^d$ and any sequence $(\nu^i)_{1\leq i\leq d}\in\{\gamma,\eta\}^d$, the random products $\prod_{i=1}^d \mathcal{W}_n^{\nu^i,N}(f_n^i)$ converge in law, as $N$ tends to infinity, to the Gaussian products $\prod_{i=1}^d \mathcal{W}_n^{\nu^i}(f_n^i)$. In addition, we have*

$$\lim_{N\to\infty}\mathbb{E}\left[\prod_{i=1}^d \mathcal{W}_n^{\nu^i,N}(f_n^i)\right] = \mathbb{E}\left[\prod_{i=1}^d \mathcal{W}_n^{\nu^i}(f_n^i)\right].$$

PROOF. We first recall from Lemmas 2.1 and 2.2 that, for $\nu\in\{\gamma,\eta\}$ and for any $f_n\in\mathcal{B}_b(F_n)$ and $p\geq 1$, we have the $\mathbb{L}_p$-mean error estimates

$$\sup_{N\geq 1}\mathbb{E}[|\mathcal{W}_n^{\nu,N}(f_n)|^p]^{1/p}<\infty$$

with the random fields $(\mathcal{W}_n^{\gamma,N},\mathcal{W}_n^{\eta,N})$ defined in (2.22) and (2.24). By the Borel–Cantelli lemma this property ensures that

$$(\gamma_n^N(f_n),\eta_n^N(f_n))\stackrel{N\to\infty}{\longrightarrow}(\gamma_n(f_n),\eta_n(f_n))\qquad\text{a.s.}$$

By the definitions of the random fields $(\mathcal{W}_n^{\gamma,N}, \mathcal{W}_n^{\gamma})$ and $(\mathcal{W}_n^{\eta,N}, \mathcal{W}_n^{\eta})$, given in (2.22), (2.23), and in (2.24), (2.25), we have that

$$\mathcal{W}_n^{\nu,N}(f_n) = \sum_{p=1}^{n} c_{p,n}^{\nu,N} \mathcal{W}_p^N(f_{p,n}^{\nu}) \quad \text{and} \quad \mathcal{W}_n^{\nu}(f_n) = \sum_{p=1}^{n} c_{p,n}^{\nu} \mathcal{W}_p(f_{p,n}^{\nu})$$

with

$$c_{p,n}^{\gamma,N} = \gamma_p^N(1) \stackrel{N \to \infty}{\longrightarrow} c_{p,n}^{\gamma} = \gamma_p(1),$$

$$c_{p,n}^{\eta,N} = \gamma_p^N(1) \; \gamma_n(1)/\gamma_n^N(1) \stackrel{N \to \infty}{\longrightarrow} c_{p,n}^{\eta} = \gamma_p(1),$$

and the pair of functions $(f_{p,n}^{\gamma}, f_{p,n}^{\eta})$ defined by

$$f_{p,n}^{\gamma} = Q_{p,n}(f_n), \qquad f_{p,n}^{\eta} = Q_{p,n}\left(\frac{1}{\gamma_n(1)}(f_n - \eta_n(f_n))\right).$$

With this notation, we find that

$$\prod_{i=1}^{d} \mathcal{W}_n^{\nu^i,N}(f_n^i) = \sum_{p_1,\ldots,p_d=1}^{n} \left[\prod_{i=1}^{d} c_{p_i,n}^{\nu^i,N}\right] \times \left[\prod_{i=1}^{d} \mathcal{W}_{p_i}^N(f_{p_i,n}^{i,\nu_i})\right],$$

$$\prod_{i=1}^{d} \mathcal{W}_n^{\nu^i}(f_n^i) = \sum_{p_1,\ldots,p_d=1}^{n} \left[\prod_{i=1}^{d} c_{p_i,n}^{\nu^i}\right] \times \left[\prod_{i=1}^{d} \mathcal{W}_{p_i}(f_{p_i,n}^{i,\nu_i})\right].$$

Recalling that the sequence of random fields $(\mathcal{W}_p^N)_{1 \leq p \leq n}$ converges in law, as $N$ tends to infinity, to a sequence of $n$ independent, Gaussian and centered random fields $(\mathcal{W}_p)_{1 \leq p \leq n}$, one concludes that $\prod_{i=1}^{d} \mathcal{W}_n^{\nu^i,N}(f_n^i)$ converges in law, as $N$ tends to infinity, to $\prod_{i=1}^{d} \mathcal{W}_n^{\nu^i}(f_n^i)$. This ends the proof of the first assertion.

Using Hölder's inequality, we can also prove that any polynomial function of terms

$$\mathcal{W}_n^{\nu,N}(f_n), \qquad \nu \in \{\gamma, \eta\}, f_n \in \mathcal{B}_b(F_n)$$

forms a uniformly integrable collection of random variables, indexed by the size and precision parameter $N \geq 1$. This property, combined with the continuous mapping theorem, and the Skorohod embedding theorem, completes the proof of the lemma. $\square$

We first present an elementary consequence of Lemma 2.4. We first rewrite (2.24) as follows

$$\mathcal{W}_n^{\eta,N}(f_n) = \mathcal{W}_n^{\gamma,N}\left(\frac{1}{\gamma_n(1)}(f_n - \eta_n(f_n))\right)$$

$$+ \left(\frac{\gamma_n(1)}{\gamma_n^N(1)} - 1\right) \times \mathcal{W}_n^{\gamma,N}\left(\frac{1}{\gamma_n(1)}(f_n - \eta_n(f_n))\right)$$

$$= \mathcal{W}_n^{\gamma,N}(\tilde{f}_n) - \frac{1}{\sqrt{N}} \frac{\gamma_n(1)}{\gamma_n^N(1)} \mathcal{W}_n^{\gamma,N}(\tilde{f}_n) \mathcal{W}_n^{\gamma,N}(\tilde{g}_n)$$

with the pair of functions $(\tilde{f}_n, \tilde{g}_n)$ defined by

$$\tilde{f}_n = \frac{1}{\gamma_n(1)}(f_n - \eta_n(f_n)) \quad \text{and} \quad \tilde{g}_n = \frac{1}{\gamma_n(1)}.$$

This yields that

$$N\mathbb{E}[\eta_n^N(f_n) - \eta_n(f_n)] = -\mathbb{E}\left[\frac{\gamma_n(1)}{\gamma_n^N(1)}\mathcal{W}_n^{\gamma,N}(\tilde{f}_n)\mathcal{W}_n^{\gamma,N}(\tilde{g}_n)\right].$$

Since the sequence of random variables $(\gamma_n(1)/\gamma_n^N(1))_{N\geq 1}$ is uniformly bounded, and it converges in law to 1, as $N$ tends to infinity, by Lemma 2.4 we conclude that

$$\begin{aligned}
(2.36) \quad \lim_{N\to\infty} N\mathbb{E}[\eta_n^N(f_n) - \eta_n(f_n)] &= -\mathbb{E}[\mathcal{W}_n^\gamma(\tilde{f}_n)\mathcal{W}_n^\gamma(\tilde{g}_n)] \\
&= -\sum_{p=1}^n \eta_p(\overline{Q}_{p,n}(1)\overline{Q}_{p,n}(f_n - \eta_n(f_n)))
\end{aligned}$$

where the renormalized semigroup $\overline{Q}_{p,n}$ is defined by

$$\overline{Q}_{p,n}(f_n) = \frac{Q_{p,n}(f_n)}{\eta_p Q_{p,n}(1)} = \frac{\gamma_p(1)}{\gamma_n(1)} Q_{p,n}(f_n).$$

We are now in position to state and prove the main result of this subsection.

THEOREM 2.5. *For any $n \geq 1$ and $\varphi_n \in \mathcal{B}_b(F_n)$, we have*

$$N\mathbb{E}[(P_n^N(a,\varphi_n) - P_n(a,\varphi_n))\mathbf{1}_{P_n^N(a)>0}]$$

$$\stackrel{N\to\infty}{\longrightarrow} -P_n(a)^{-2} \sum_{p=1}^n \gamma_p(1)\gamma_p^-(1)\eta_p^-[P_{p,n}(a)P_{p,n}(a,\varphi_n)]$$

*with the collection of functions $P_{p,n}(a)$, $P_{p,n}(a,\varphi_n) \in \mathcal{B}_b(F_p)$ defined, respectively, in* (2.28) *and* (2.31).

PROOF. The proof is essentially based on a judicious way to rewrite (2.34). If we define

$$f_n^{(a)} = \frac{T_n^{(a)}}{\eta_n T_n^{(a)}(1)}(\varphi_n - P_n(a,\varphi_n)) \quad \text{and} \quad g_n^{(a)} = \frac{T_n^{(a)}(1)}{\eta_n T_n^{(a)}(1)},$$

then, on the event $\{P_n^N(a) > 0\}$, we have

$$N[P_n^N(a, \varphi_n) - P_n(a, \varphi_n)]$$
$$= N[\eta_n^N(f_n^{(a)}) - \eta_n(f_n^{(a)})] - \frac{1}{\eta_n^N(g_n^{(a)})} \mathcal{W}_n^{\eta,N}(f_n^{(a)}) \mathcal{W}_n^{\eta,N}(g_n^{(a)}).$$

By Lemma 2.4 and (2.36) we conclude that

$$N\mathbb{E}[(P_n^N(a, \varphi_n) - P_n(a, \varphi_n))\mathbf{1}_{P_n^N(a)>0}]$$
$$\stackrel{N\to\infty}{\longrightarrow} -\mathbb{E}[\mathcal{W}_n^\eta(f_n^{(a)})\mathcal{W}_n^\eta(g_n^{(a)})] - \mathbb{E}\left[\mathcal{W}_n^\gamma\left(\frac{f_n^{(a)}}{\gamma_n(1)}\right)\mathcal{W}_n^\gamma\left(\frac{1}{\gamma_n(1)}\right)\right].$$

On the other hand, using (2.25) we find that

$$\mathbb{E}[\mathcal{W}_n^\eta(f_n^{(a)})\mathcal{W}_n^\eta(g_n^{(a)})]$$
$$= \sum_{p=1}^n (\gamma_p(1)/\gamma_n(1))^2 \mathbb{E}[\mathcal{W}_p(Q_{p,n}(f_n^{(a)}))\mathcal{W}_p(Q_{p,n}(g_n^{(a)} - 1))]$$
$$= \sum_{p=1}^n (\gamma_p(1)/\gamma_n(1))^2 \eta_p(Q_{p,n}(f_n^{(a)})Q_{p,n}(g_n^{(a)} - 1)).$$

Similarly, by (2.23) we have

$$\mathbb{E}\left[\mathcal{W}_n^\gamma\left(\frac{f_n^{(a)}}{\gamma_n(1)}\right)\mathcal{W}_n^\gamma\left(\frac{1}{\gamma_n(1)}\right)\right]$$
$$= \sum_{p=1}^n \gamma_p(1)^2 \, \mathbb{E}\left[\mathcal{W}_p\left(Q_{p,n}\frac{f_n^{(a)}}{\gamma_n(1)}\right)\mathcal{W}_p\left(Q_{p,n}\frac{1}{\gamma_n(1)}\right)\right].$$

It is now convenient to notice that

$$\mathbb{E}\left[\mathcal{W}_p\left(Q_{p,n}\frac{f_n^{(a)}}{\gamma_n(1)}\right)\mathcal{W}_p\left(Q_{p,n}\frac{1}{\gamma_n(1)}\right)\right]$$
$$= \gamma_n(1)^{-2}\mathbb{E}[\mathcal{W}_p(Q_{p,n}(f_n^{(a)}))\mathcal{W}_p(Q_{p,n}(1))]$$
$$= \gamma_n(1)^{-2} \times \eta_p(Q_{p,n}(f_n^{(a)})Q_{p,n}(1)).$$

This implies that

$$\mathbb{E}\left[\mathcal{W}_n^\gamma\left(\frac{f_n^{(a)}}{\gamma_n(1)}\right)\mathcal{W}_n^\gamma\left(\frac{1}{\gamma_n(1)}\right)\right] = \sum_{p=1}^n (\gamma_p(1)/\gamma_n(1))^2 \eta_p(Q_{p,n}(1)Q_{p,n}(f_n^{(a)}))$$

from which we conclude that

(2.37)
$$N\mathbb{E}([P_n^N(a, \varphi_n) - P_n(a, \varphi_n)]\mathbf{1}_{P_n^N(a)>0})$$
$$\stackrel{N\to\infty}{\longrightarrow} -\sum_{p=1}^n (\gamma_p(1)/\gamma_n(1))^2 \eta_p(Q_{p,n}(f_n^{(a)})Q_{p,n}(g_n^{(a)})).$$

By the definition of the function $T_n^{(a)}(\varphi_n)$ we have $\eta_n T_n^{(a)}(1) = P_n(a)/\gamma_n(1)$ and for any $y_p = (x_0, \ldots, x_p) \in F_p$

$$Q_{p,n}(T_n^{(a)}(\varphi_n))(x_0, \ldots, x_p)$$
$$= \left[\prod_{1 \leq k < p} G_k^-(x_0, \ldots, x_k)\right]$$
$$\times \mathbb{E}[\varphi_n(X_0, \ldots, X_n)\mathbf{1}_{V_n(X_n) \geq a}|(X_0, \ldots, X_p) = (x_0, \ldots, x_p)].$$

By the definition of the pair of functions $(f_n^{(a)}, g_n^{(a)})$, these observations yield

$$Q_{p,n}(f_n^{(a)})(x_0, \ldots, x_p) = \frac{\gamma_n(1)}{P_n(a)}\left[\prod_{1 \leq k < p} G_k^-(x_0, \ldots, x_k)\right]$$
$$\times P_{p,n}(a, \varphi_n)(x_0, \ldots, x_p),$$

$$Q_{p,n}(g_n^{(a)})(x_0, \ldots, x_p) = \frac{\gamma_n(1)}{P_n(a)}\left[\prod_{1 \leq k < p} G_k^-(x_0, \ldots, x_k)\right]P_{p,n}(a)(x_0, \ldots, x_p).$$

To take the final step, we notice that

$$Q_{p,n}(f_n^{(a)})(x_0, \ldots, x_p) \times Q_{p,n}(g_n^{(a)})(x_0, \ldots, x_p) \times (P_n(a)/\gamma_n(1))^2$$
$$= \left[\prod_{1 \leq k < p} G_k^-(x_0, \ldots, x_k)^2\right]P_{p,n}(a, \varphi_n)(x_0, \ldots, x_p)P_{p,n}(a)(x_0, \ldots, x_p).$$

This implies that

$$\gamma_p(Q_{p,n}(f_n^{(a)})Q_{p,n}(g_n^{(a)})) \times (P_n(a)/\gamma_n(1))^2$$
$$= \mathbb{E}\left(\left[\prod_{1 \leq k < p} G_k^-(X_0, \ldots, X_k)\right]\right.$$
$$\left.\times P_{p,n}(a, \varphi_n)(X_0, \ldots, X_p)P_{p,n}(a)(X_0, \ldots, X_p)\right)$$
$$= \gamma_p^-(P_{p,n}(a, \varphi_n)P_{p,n}(a)) = \gamma_p^-(1) \times \eta_p^-(P_{p,n}(a, \varphi_n)P_{p,n}(a)).$$

Using the identity $\gamma_p = \gamma_p(1)\eta_p$ and substituting the last equation into (2.37) completes the proof. $\square$

2.7. *Variance comparisons for Gaussian particle models.* Let $(X_p)_{1 \leq p \leq n}$ be the Gaussian sequence defined in (2.2). We consider the elementary energy-like function $V_n(x) = x$, and the Feynman–Kac twisted models associated to the potential functions

(2.38) $\quad G_p(x_0, \ldots, x_p) = \exp[\lambda(x_p - x_{p-1})] \quad$ for some $\lambda > 0$.

Arguing as in (2.4), we prove that the Feynman–Kac distribution $\eta_p^-$ is the path distribution of the chain defined by the recursion

$$X_p^- = X_{p-1}^- + W_p \quad \text{and} \quad X_k^- = X_{k-1}^- - \lambda + W_k, \qquad 1 \leq k < p,$$

where $X_0 = 0$, and where $(W_k)_{1 \leq k \leq p}$ represents a sequence of independent and identically distributed Gaussian random variables, with $(\mathbb{E}(W_1), \mathbb{E}(W_1^2)) = (0, 1)$. We also observe that in this case we have

$$(2.39) \qquad \gamma_p(1)\gamma_p^-(1) = \mathbb{E}[e^{\lambda X_{p-1}}]^2 = e^{\lambda^2(p-1)}.$$

The next lemma is instrumental for estimating the quantities $\eta_p^-(P_{p,n}(a)^2)$ introduced in (2.27).

LEMMA 2.6. *Let $(W_1, W_2)$ be a pair of independent Gaussian random variables, with $(\mathbb{E}(W_i), \mathbb{E}(W_i^2)) = (0, \sigma_i^2)$, with $\sigma_i > 0$ and $i = 1, 2$. Then, for any $a > 0$, we have the exponential estimate*

$$C(a, \sigma_1, \sigma_2) \leq \mathbb{E}[\mathbb{P}(W_1 + W_2 \geq a | W_1)^2] \exp\left(\frac{a^2}{2\sigma_1^2 + \sigma_2^2}\right) \leq 1,$$

*where*

$$C(a, \sigma_1, \sigma_2) = (2\pi)^{3/2} \left(\frac{\sigma_2 a}{2\sigma_1^2 + \sigma_2^2} + \frac{2\sigma_1^2 + \sigma_2^2}{\sigma_2 a}\right)^{-2} \left(\frac{2\sigma_1 a}{2\sigma_1^2 + \sigma_2^2} + \frac{2\sigma_1^2 + \sigma_2^2}{2\sigma_1 a}\right)^{-1}.$$

PROOF. Using exponential version of Chebyshev's inequality we first check that, for any $\lambda > 0$, we have

$$\mathbb{P}(W_1 + W_2 \geq a | W_1) \leq e^{\lambda(W_1 - a)} \mathbb{E}(e^{\lambda W_2}) = e^{\lambda(W_1 - a) + \lambda^2 \sigma_2^2/2}.$$

Integrating the random variable $W_1$ and choosing $\lambda = a/(2\sigma_1^2 + \sigma_2^2)$, we establish the upper bound

$$\mathbb{E}[\mathbb{P}(W_1 + W_2 \geq a | W_1)^2] \leq e^{-2\lambda a + \lambda^2(2\sigma_1^2 + \sigma_2^2)} = e^{-a^2/(2\sigma_1^2 + \sigma_2^2)}.$$

For any $\varepsilon \in (0, 1)$, we have

$$\mathbb{E}[\mathbb{P}(W_1 + W_2 \geq a | W_1)^2] \geq \mathbb{P}(W_2 \geq \varepsilon a)^2 \mathbb{P}(W_1 \geq (1 - \varepsilon)a).$$

Applying Mill's inequality yields

$$\mathbb{E}[\mathbb{P}(W_1 + W_2 \geq a | W_1)^2] \geq \frac{(2\pi)^{3/2}}{(\varepsilon a/\sigma_2 + \sigma_2/(\varepsilon a))^2((1-\varepsilon)a/\sigma_1 + \sigma_1/((1-\varepsilon)a))} \\ \times e^{-a^2(\varepsilon^2/\sigma_2^2 + (1-\varepsilon)^2/(2\sigma_1^2))}.$$

Choosing $\varepsilon = \sigma_2^2/(2\sigma_1^2 + \sigma_2^2)$ establishes the lower bound. □

From previous considerations, we notice that
$$\eta_p^-(P_{p,n}(a)^2) = \mathbb{E}[\mathbb{P}(W_1 + W_2 \geq (a + \lambda(p-1))|W_1)^2],$$
where $(W_1, W_2)$ are a pair of independent and centered Gaussian random variables, with $(\mathbb{E}(W_1^2), \mathbb{E}(W_2^2)) = (p, n-p)$. Lemma 2.6 now implies that
$$(2.40) \qquad \eta_p^-(P_{p,n}(a)^2) \leq \exp[-(a + \lambda(p-1))^2/(n+p)].$$
Substituting the estimates (2.39) and (2.40) into (2.27), we find that
$$\sigma_n^\gamma(a)^2 \leq \sum_{p=1}^n [e^{\lambda^2(p-1) - (a+\lambda(p-1))^2/(n+p)} - P_n(a)^2]$$
$$= \sum_{0 \leq p < n} [e^{-a^2/n} e^{(p+1)/(n(n+p+1))[a - \lambda(np)/(p+1)]^2 + \lambda^2 p/(p+1)} - P_n(a)^2].$$
For $\lambda = a/n$, this yields that
$$\sigma_n^\gamma(a)^2 \leq \sum_{0 \leq p < n} [e^{-a^2/n} e^{a^2/n^2(1 - 1/(n+p+1))} - P_n(a)^2]$$
$$(2.41) \qquad \leq n(e^{-(a^2/n)(1-1/n)} - P_n(a)^2).$$

We find that this estimate has the same exponential decay as the one obtained in (2.6) for the corresponding noninteracting IS model. The only difference between these two asymptotic variances comes from the multiplication parameter $n$. This additional term can be interpreted as the number of interactions used in the construction of the genealogical tree simulation model. We can compare the efficiencies of the IPS strategy and the usual MC strategy which are two methods that do not require to twist the input probability distribution, in contrast to IS. The IPS provides a variance reduction by a factor of the order of $nP_n(a)$. In practice the number $n$ of selection steps is of the order of ten or a few tens, while the goal is the estimation of a probability $P_n(a)$ of the order of $10^{-6}$–$10^{-12}$. The gain is thus very significant, as we shall see in the numerical applications.

Now, we consider the Feynman–Kac twisted models associated to the potential functions
$$(2.42) \qquad G_p(x_0, \ldots, x_p) = \exp(\lambda x_p) \qquad \text{for some } \lambda > 0.$$
Arguing as in (2.11), we prove that $\eta_p^-$ is the distribution of the Markov chain
$$X_p^- = X_{p-1}^- + W_p \quad \text{and} \quad X_k^- = X_{k-1}^- - \lambda(p-k) + W_k, \qquad 1 \leq k < p,$$
where $X_0 = 0$, and where $(W_k)_{1 \leq k \leq p}$ represents a sequence of independent and identically distributed Gaussian random variables, with $(\mathbb{E}(W_1), \mathbb{E}(W_1^2)) = (0, 1)$. We also notice that
$$(2.43) \qquad \gamma_p(1)\gamma_p^-(1) = \mathbb{E}[e^{\lambda \sum_{1 \leq k < p} X_k}]^2 = e^{\lambda^2 \sum_{1 \leq k < p} k^2}.$$

In this situation, we observe that

$$\eta_p^-(P_{p,n}(a)^2) = \mathbb{E}\left[\mathbb{P}\left(W_1 + W_2 \geq a + \lambda \sum_{1 \leq k < p} k | W_1\right)^2\right],$$

where $(W_1, W_2)$ are a pair of independent and centered Gaussian random variables, with $(\mathbb{E}(W_1^2), \mathbb{E}(W_2^2)) = (p, n-p)$. As before, Lemma 2.6 now implies that

(2.44) $$\eta_p^-(P_{p,n}(a)^2) \leq \exp\left[-\frac{1}{n+p}\left(a + \lambda \frac{p(p-1)}{2}\right)^2\right].$$

Using the estimates (2.43) and (2.44), and recalling that $\sum_{1 \leq k \leq n} k^2 = n(n+1)(2n+1)/6$, we conclude that

$$\sigma_n^\gamma(a)^2 \leq \sum_{p=1}^n [e^{(1/6)\lambda^2(p-1)p(2p-1) - (a+\lambda p(p-1)/2)^2/(n+p)} - P_n(a)^2]$$

$$= \sum_{p=1}^n [e^{-(a^2/n)} e^{p/(n(n+p))[a - \lambda n(p-1)/2]^2 + (1/12)\lambda^2(p-1)p(p+1)} - P_n(a)^2].$$

If we take $\lambda = 2a/[n(n-1)]$, then we get

$$\sigma_n^\gamma(a)^2 \leq \sum_{p=1}^n [e^{-a^2/n} e^{a^2/(n^2(n-1)^2)[np/(n+p)(n-p)^2 + (p-1)p(p+1)/3]} - P_n(a)^2]$$

$$= \sum_{p=1}^n [e^{-a^2/n} e^{(a^2/n)n^2/(n-1)^2 [\theta(p/n) - p/(3n^3)]} - P_n(a)^2]$$

with the increasing function $\theta : \varepsilon \in [0,1] \longrightarrow \theta(\varepsilon) = \varepsilon\frac{(1-\varepsilon)^2}{(1+\varepsilon)} + \frac{\varepsilon^3}{3} \in [0, 1/3]$. From these observations, we deduce the estimate

(2.45) $$\sigma_n^\gamma(a)^2 \leq n[e^{-(a^2/n)(2/3)(1-1/(n-1))} - P_n(a)^2].$$

Note that the inequalities are sharp in the exponential sense by the lower bound obtained in Lemma 2.6. Accordingly we get that the asymptotic variance is not of the order of $P_n(a)^2$, but rather $P_n(a)^{4/3}$. As in the first Gaussian example, we observe that this estimate has the same exponential decays as the one obtained in (2.12) for the corresponding IS algorithm. But, once again, the advantage of the IPS method compared to IS is that it does not require to twist the original transition probabilities, which makes the IPS strategy much easier to implement.

*Conclusion.* The comparison of the variances (2.41) and (2.45) shows that the variance of the estimator $P_n^N(a)$ is much smaller when the potential (2.38) is used rather than the potential (2.42). We thus get the important qualitative conclusion that it is not efficient to select the "best" particles (i.e., those with the highest energy values), but it is much more efficient to select amongst the particles with the best energy increments. This conclusion is also an a posteriori justification of the necessity to carry out a path-space analysis, and not only a state-space analysis. The latter one is simpler but it cannot consider potentials of the form (2.38) that turn out to be more efficient.

**3. Estimation of the tail of a probability density function.** We collect and sum up the general results presented in Section 2 and we apply them to propose an estimator for the tail of the probability density function (p.d.f.) of a real-valued function of a Markov chain. We consider an $(E, \mathcal{E})$-valued Markov chain $(X_p)_{0 \le p \le n}$ with nonhomogeneous transition kernels $K_p$. In a first time, we show how the results obtained in the previous section allow us to estimate the probability of a rare event of the form $\{V(X_n) \in A\}$:

$$(3.1) \qquad P_A = \mathbb{P}(V(X_n) \in A) = \mathbb{E}[\mathbf{1}_A(V(X_n))],$$

where $V$ is some function from $E$ to $\mathbb{R}$. We shall construct an estimator based on an IPS. As pointed out in the previous section, the quality of the estimator depends on the choice on the weight function. The weight function should fulfill two conditions. First, it should favor the occurrence of the rare event without involving too large normalizing constants. Second, it should give rise to an algorithm that can be easily implemented. Indeed the implementation of the IPS with an arbitrary weight function requires recording the complete set of path-particles. If $N$ particles are generated and time runs from 0 to $n$, this set has size $(n+1) \times N \times \dim(E)$ which may exceed the memory capacity of the computer. The weight function should be chosen so that only a smaller set needs to be recorded to compute the estimator of the probability of occurrence of the rare event. We shall examine two weight functions and the two corresponding algorithms that fulfill both conditions.

ALGORITHM 1. Let us fix some $\beta > 0$. The first algorithm is built with the weight function

$$(3.2) \qquad G_p^\beta(x) = \exp[\beta V(x_p)].$$

The practical implementation of the IPS reads as follows.

*Initialization.* We start with a set of $N$ i.i.d. initial conditions $\hat{X}_0^{(i)}$, $1 \leq i \leq N$, chosen according to the initial distribution of $X_0$. This set is complemented with a set of weights $\hat{Y}_0^{(i)} = 1$, $1 \leq i \leq N$. This forms a set of $N$ particles: $(\hat{X}_0^{(i)}, \hat{Y}_0^{(i)})$, $1 \leq i \leq N$, where a particle is a pair $(\hat{X}, \hat{Y}) \in E \times \mathbb{R}^+$.

Now, assume that we have a set of $N$ particles at time $p$ denoted by $(\hat{X}_p^{(i)}, \hat{Y}_p^{(i)})$, $1 \leq i \leq N$.

*Selection.* We first compute the normalizing constant

$$(3.3) \qquad \hat{\eta}_p^N = \frac{1}{N} \sum_{i=1}^N \exp[\beta V(\hat{X}_p^{(i)})].$$

We choose independently $N$ particles according to the empirical distribution

$$(3.4) \qquad \mu_p^N(d\check{X}, d\check{Y}) = \frac{1}{N\hat{\eta}_p^N} \times \sum_{i=1}^N \exp[\beta V(\hat{X}_p^{(i)})] \delta_{(\hat{X}_p^{(i)}, \hat{Y}_p^{(i)})}(d\check{X}, d\check{Y}).$$

The new particles are denoted by $(\check{X}_p^{(i)}, \check{Y}_p^{(i)})$, $1 \leq i \leq N$.

*Mutation.* For every $1 \leq i \leq N$, the particle $(\check{X}_p^{(i)}, \check{Y}_p^{(i)})$ is transformed into $(\hat{X}_{p+1}^{(i)}, \hat{Y}_{p+1}^{(i)})$ by the mutation procedure

$$(3.5) \qquad \check{X}_p^{(i)} \xrightarrow{K_{p+1}} \hat{X}_{p+1}^{(i)},$$

where the mutations are performed independently, and

$$(3.6) \qquad \hat{Y}_{p+1}^{(i)} = \check{Y}_p^{(i)} \exp[-\beta V(\check{X}_p^{(i)})].$$

The memory required by the algorithm is $N \dim(E) + N + n$, where $N \dim(E)$ is the memory required by the record of the set of particles, $N$ is the memory required by the record of the set of weights and $n$ is the memory required by the record of the normalizing constants $\hat{\eta}_p^N$, $0 \leq p \leq n-1$. The estimator of the probability $P_A$ is then

$$(3.7) \qquad P_A^N = \left[\frac{1}{N} \sum_{i=1}^N \mathbf{1}_A(V(\hat{X}_n^{(i)})) \hat{Y}_n^{(i)}\right] \times \prod_{k=0}^{n-1} \hat{\eta}_p^N.$$

This estimator is unbiased in the sense that $\mathbb{E}[P_A^N] = P_A$. The central limit theorem for the estimator states that

$$(3.8) \qquad \sqrt{N}(P_A^N - P_A) \xrightarrow{N \to \infty} \mathcal{N}(0, Q_A)$$

where the variance is

$$Q_A = \sum_{p=1}^n \mathbb{E}\left[\mathbb{E}_{X_p}[\mathbf{1}_A(V(X_n))]^2 \prod_{k=0}^{p-1} G_k^{-1}(X)\right] \mathbb{E}\left[\prod_{k=0}^{p-1} G_k(X)\right]$$
$$(3.9) \quad - \mathbb{E}[\mathbf{1}_A(X_n)]^2.$$

ALGORITHM 2. Let us fix some $\alpha > 0$. The second algorithm is built with the weight function

$$(3.10) \quad G_p^\alpha(x) = \exp[\alpha(V(x_p) - V(x_{p-1}))].$$

*Initialization.* We start with a set of $N$ i.i.d. initial conditions $\hat{X}_0^{(i)}$, $1 \leq i \leq N$, chosen according to the initial distribution of $X_0$. This set is complemented with a set of parents $\hat{W}_0^{(i)} = x_0$, $1 \leq i \leq N$, where $x_0$ is an arbitrary point of $E$ with $V(x_0) = V_0$. This forms a set of $N$ particles: $(\hat{W}_0^{(i)}, \hat{X}_0^{(i)})$, $1 \leq i \leq N$, where a particle is a pair $(\hat{W}, \hat{X}) \in E \times E$.

Now, assume that we have a set of $N$ particles at time $p$ denoted by $(\hat{W}_p^{(i)}, \hat{X}_p^{(i)})$, $1 \leq i \leq N$.

*Selection.* We first compute the normalizing constant

$$(3.11) \quad \hat{\eta}_p^N = \frac{1}{N} \sum_{i=1}^N \exp[\alpha(V(\hat{X}_p^{(i)}) - V(\hat{W}_p^{(i)}))].$$

We choose independently $N$ particles according to the empirical distribution

$$\mu_p^N(d\check{W}, d\check{X}) = \frac{1}{N\hat{\eta}_p^N} \sum_{i=1}^N \exp[\alpha(V(\hat{X}_p^{(i)}) - V(\hat{W}_p^{(i)}))]$$
$$(3.12) \quad \times \delta_{(\hat{W}_p^{(i)}, \hat{X}_p^{(i)})}(d\check{W}, d\check{X}).$$

The new particles are denoted by $(\check{W}_p^{(i)}, \check{X}_p^{(i)})$, $1 \leq i \leq N$.

*Mutation.* For every $1 \leq i \leq N$, the particle $(\check{W}_p^{(i)}, \check{X}_p^{(i)})$ is transformed into $(\hat{W}_{p+1}^{(i)}, \hat{X}_{p+1}^{(i)})$ by the mutation procedure $\check{X}_p^{(i)} \stackrel{K_{p+1}}{\longrightarrow} \hat{X}_{p+1}^{(i)}$ where the mutations are performed independently, and $\hat{W}_{p+1}^{(i)} = \check{X}_p^{(i)}$.

The memory required by the algorithm is $2N \dim(E) + n$. The estimator of the probability $P_A$ is then

$$(3.13) \quad P_A^N = \left[\frac{1}{N} \sum_{i=1}^N \mathbf{1}_A(V(\hat{X}_n^{(i)})) \exp(-\alpha(V(\hat{W}_n^{(i)}) - V_0))\right] \times \left[\prod_{k=0}^{n-1} \hat{\eta}_p^N\right].$$

This estimator is unbiased and satisfies the central limit theorem (3.8).

Let us now focus our attention to the estimation of the p.d.f. tail of $V(X_n)$. The rare event is then of the form $\{V(X_n) \in [a, a+\delta a)\}$ with a large $a$ and an evanescent $\delta a$. We assume that the p.d.f. of $V(X_n)$ is continuous so that the p.d.f. can be seen as

$$p(a) = \lim_{\delta a \to 0} \frac{1}{\delta a} p_{\delta a}(a), \qquad p_{\delta a}(a) = \mathbb{P}(V(X_n) \in [a, a+\delta a)).$$

We propose to use the estimator

$$(3.14) \qquad p_{\delta a}^N(a) = \frac{1}{\delta a} \times P_{[a, a+\delta a)}^N$$

with a small $\delta a$. The central limit theorem for the p.d.f. estimator takes the form

$$(3.15) \qquad \sqrt{N}(p_{\delta a}^N(a) - p_{\delta a}(a)) \stackrel{N \to \infty}{\longrightarrow} \mathcal{N}(0, p_2^2(a, \delta a)),$$

where the variance $p_2^2(a, \delta a)$ has a limit $p_2^2(a)$ as $\delta a$ goes to 0 which admits a simple representation formula:

$$(3.16) \quad p_2^2(a) = \lim_{\delta a \to 0} \frac{1}{\delta a} \mathbb{E}\left[\mathbf{1}_{[a,a+\delta a)}(V(X_n)) \prod_{k=0}^{n-1} G_k^{-1}(X)\right] \mathbb{E}\left[\prod_{k=0}^{n-1} G_k(X)\right].$$

Note that all other terms in the sum (3.9) are of order $\delta a^2$ and are therefore negligible. This is true as soon as the distribution of $V(X_n)$ given $X_p$ for $p < n$ admits a bounded density with respect to the Lebesgue measure. Accordingly, the variance $p_2^2(a)$ can be estimated by $\lim_{\delta a \to 0}(\delta a)^{-1} Q_{[a,a+\delta a)}^N$, where $Q_A^N$ is given by

$$(3.17) \qquad Q_A^N = \left[\frac{1}{N}\sum_{i=1}^N \mathbf{1}_A(V(\hat{X}_n^{(i)}))(\hat{Y}_n^{(i)})^2\right] \times \left[\prod_{k=0}^{n-1} \hat{\eta}_p^N\right]^2$$

for Algorithm 1, and by

$$(3.18) \quad Q_A^N = \left[\frac{1}{N}\sum_{i=1}^N \mathbf{1}_A(V(\hat{X}_n^{(i)})) \exp(-2\alpha(V(\hat{W}_n^{(i)}) - V_0))\right] \times \left[\prod_{k=0}^{n-1} \hat{\eta}_p^N\right]^2$$

for Algorithm 2. The estimators of the variances are important because confidence intervals can then be obtained.

The variance analysis carried out in Section 2.7 predicts that the second algorithm [with the potential (3.2)] should give better results than the Algorithm 1 [with the potential (3.10)]. We are going to illustrate this important statement in the following sections devoted to numerical simulations.

From a practical point of view, it can be interesting to carry out several IPS simulations with different values for the parameters $\beta$ (Algorithm 1) and $\alpha$ (Algorithm 2), and also one MC simulation. It is then possible to reconstruct the p.d.f. of $V(X_n)$ by the following procedure. Each IPS or MC simulation gives an estimation for the p.d.f. $p$ (whose theoretical value does not depend on the method) and also an estimation for the ratio $p_2/p$ (whose theoretical value depends on the method). We first consider the different estimates of $a \mapsto p_2/p(a)$ and detect, for each given value of $a$, which IPS gives the minimal value of $p_2/p(a)$. For this value of $a$, we then use the estimation of $p(a)$ obtained with this IPS. This method will be used in Section 5.

**4. A toy model.** In this section we apply the IPS method to compute the probabilities of rare events for a very simple system for which we know explicit formulas. The system under consideration is the Gaussian random walk $X_{p+1} = X_p + W_{p+1}$, $X_0 = 0$, where the $(W_p)_{p=1,\ldots,n}$ are i.i.d. Gaussian random variables with zero mean and variance 1. Let $n$ be some positive integer. The goal is to compute the p.d.f. of $X_n$, and in particular the tail corresponding to large positive values.

We choose the weight function

$$G_p^\alpha(x) = \exp[\alpha(x_p - x_{p-1})]. \tag{4.1}$$

The theoretical p.d.f. is a Gaussian p.d.f. with variance $n$:

$$p(a) = \frac{1}{\sqrt{2\pi n}} \exp\left(-\frac{a^2}{2n}\right). \tag{4.2}$$

The theoretical variance of the p.d.f. estimator can be computed from (3.16):

$$p_2^2(a) = p^2(a) \times \sqrt{2\pi n} \exp\left(\alpha^2 \frac{n-1}{n} + \frac{(a - \alpha(n-1))^2}{2n}\right). \tag{4.3}$$

When $\alpha = 0$, we have $p_2^2(a) = p(a)$, which is the result of standard MC. For $\alpha \neq 0$, the ratio $p_2(a)/p(a)$ is minimal when $a = \alpha(n-1)$ and then $p_2(a) = p(a) \sqrt[4]{2\pi n} \exp(\alpha^2(n-1)/(2n))$. This means that the IPS with some given $\alpha$ is especially relevant for estimating the p.d.f. tail around $a = \alpha(n-1)$.

Let us assume that $n \gg 1$. Typically we look for the p.d.f. tail for $a = a_0 \sqrt{n}$ with $a_0 > 1$ because $\sqrt{n}$ is the typical value of $X_n$. The optimal choice is $\alpha = a_0/\sqrt{n}$ and then the relative error is $p_2(a)/p(a) \simeq \sqrt[4]{2\pi n}$.

In Figure 1 we compare the results from MC simulations, IPS simulations and theoretical formulas with the weight function (4.1). We use a set of $2 \times 10^4$ particles to estimate the p.d.f. tail of $X_n$ with $n = 15$. The agreement shows that we can be confident with the results given by the IPS for predicting rare events with probabilities $10^{-12}$.

We now choose the weight function

$$G_p^\beta(x) = \exp(\beta x_p). \tag{4.4}$$

We get the same results, but the explicit expression for the theoretical variance of the p.d.f. estimator is

$$(4.5) \quad p_2^2(a) = p^2(a) \times \sqrt{2\pi n} \exp\left(\beta^2 \frac{n(n^2-1)}{12} + \frac{(a - \beta n(n-1)/2)^2}{2n}\right).$$

When $\beta = 0$, we have $p_2^2(a) = p(a)$, which is the result of standard MC. For $\beta \neq 0$, the ratio $p_2(a)/p(a)$ is minimal when $a = \beta n(n-1)/2$ and then $p_2(a) = p(a)\sqrt[4]{2\pi n} \exp(\beta^2 n(n^2-1)/24)$. This means that the IPS with some given $\beta$ is especially relevant for estimating the p.d.f. tail around $a = \beta n(n-1)/2$.

Let us assume that $n \gg 1$. Typically we look for the p.d.f. tail for $a = a_0\sqrt{n}$ with $a_0 > 1$. The optimal choice is $\beta = 2a_0/n^{3/2}$ and then the relative error is $p_2(a)/p(a) \simeq (2\pi n)^{1/4} \exp(a_0^2/6) = (2\pi n)^{-1/12} p(a)^{-1/3}$. The relative error is larger than the one we get with the weight function (4.1). In Figure 2 we compare the results from MC simulations, IPS simulations and the theoretical formulas with the weight function (4.4). This shows that the weight function (4.4) is less efficient than (4.1). Thus the numerical simulations confirm the variance comparison carried out in Section 2.7.

## 5. Polarization mode dispersion in optical fibers.

5.1. *Introduction.* The study of pulse propagation in a fiber with random birefringence has become of great interest for telecommunication applications. Recent experiments have shown that polarization mode dispersion (PMD) is one of the main limitations on fiber transmission links because it can involve significant pulse broadening [13]. PMD has its origin in the birefringence [27], that is, the fact that the electric field is a vector field and the index of refraction of the medium depends on the polarization state (i.e., the unit vector pointing in the direction of the electric vector field). Random birefringence results from variations of the fiber parameters such as the core radius or geometry. There exist various physical reasons for the fluctuations of the fiber parameters. They may be induced by mechanical distortions on fibers in practical use, such as pointlike pressures or twists [21]. They may also result from variations of ambient temperature or other environmental parameters [2].

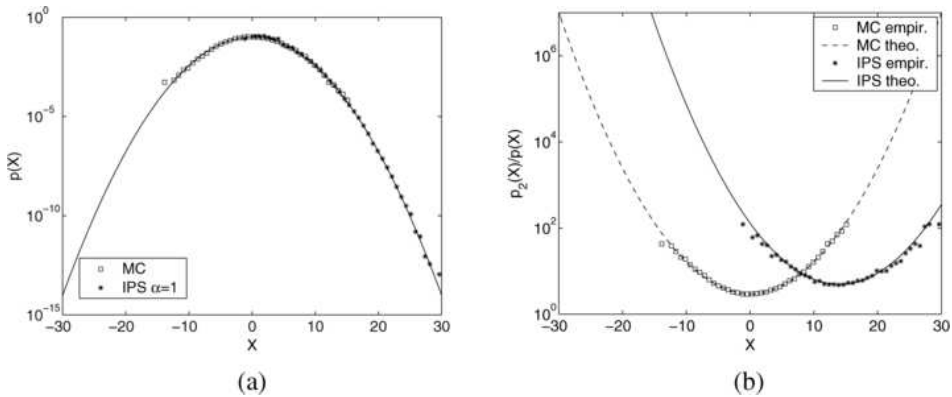

FIG. 1. (a) *P.d.f. estimations obtained by the usual MC technique (squares) and by the IPS with the weight function* (4.1) *with* $\alpha = 1$ *(stars). The solid line stands for the theoretical Gaussian distribution.* (b) *Empirical and theoretical ratios* $p_2/p$.

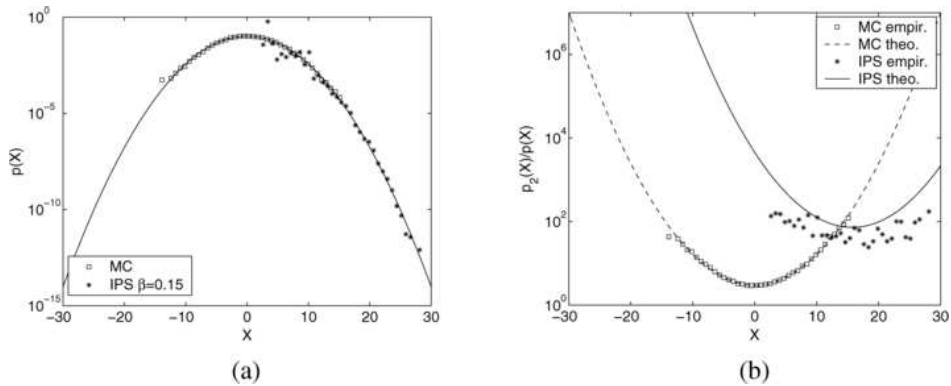

FIG. 2. (a) *P.d.f. estimations obtained by the usual MC technique (squares) and by the IPS with the weight function* (4.4) *with* $\beta = 0.15$ *(stars). The solid line stands for the theoretical Gaussian distribution.* (b) *Empirical and theoretical ratios* $p_2/p$.

The difficulty is that PMD is a random phenomenon. Designers want to ensure that some exceptional but very annoying event occurs only a very small fraction of time. This critical event corresponds to a pulse spreading beyond a threshold value. For example, a designer might require that such an event occurs less than 1 minute per year [3]. PMD in an optical fiber varies with time due to vibrations and variations of environmental parameters. The usual assumption is that the fiber passes ergodically through all possible realizations. Accordingly requiring that an event occurs a fraction of time $p$ is equivalent to requiring that the probability of this event is $p$. The problem is then reduced to the estimation of the probability of a rare event. Typically the probability is $10^{-6}$ or less [3]. It is extremely difficult to use either laboratory experiments or MC simulations to obtain a reliable estimate of such a low probability because the number of configurations that must be explored is very large. Recently IS has been applied to numerical simulations of PMD [2]. This method gives good results; however, it requires very good physical insight into the problem because it is necessary for the user to know how to produce artificially large pulse widths. We would like to revisit this work by applying the IPS strategy. The main advantage is that we do not need to specify how to produce artificially large pulse widths, as the IPS will automatically select the good "particles."

5.2. *Review of PMD models.* The pulse spreading in a randomly birefringent fiber is characterized by the so-called square differential group delay (DGD) $\tau = |\hat{\mathbf{r}}|^2$. The vector $\hat{\mathbf{r}}$ is the so-called PMD vector, which is solution of

$$\hat{\mathbf{r}}_z = \omega \Omega(z) \times \hat{\mathbf{r}} + \Omega(z), \tag{5.1}$$

where $\Omega(z)$ is a three-dimensional zero-mean stationary random process modeling PMD.

5.2.1. *The white noise model.* Simplified analytical models have been studied. In the standard model [12, 13, 20, 27] it is assumed that the process $\Omega$ is a white noise with autocorrelation function $\mathbb{E}[\Omega_i(z')\Omega_j(z)] = \sigma^2 \delta_{ij} \delta(z' - z)$. In such a case the differential equation (5.1) must be interpreted as a stochastic differential equation:

$$d\hat{r}_1 = \sigma\omega \hat{r}_3 \circ dW_z^2 - \sigma\omega \hat{r}_2 \circ dW_z^3 + \sigma\, dW_z^1, \tag{5.2}$$

$$d\hat{r}_2 = \sigma\omega \hat{r}_1 \circ dW_z^3 - \sigma\omega \hat{r}_3 \circ dW_z^1 + \sigma\, dW_z^2, \tag{5.3}$$

$$d\hat{r}_3 = \sigma\omega \hat{r}_2 \circ dW_z^1 - \sigma\omega \hat{r}_1 \circ dW_z^2 + \sigma\, dW_z^3, \tag{5.4}$$

where $\circ$ stands for the Stratonovich integral and the $W^j$'s are three independent Brownian motions. It is then possible to establish [12] that the DGD $\tau$ is a diffusion process and in particular that $\tau(\omega, z)$ obeys a Maxwellian distribution if $\hat{\mathbf{r}}(0) = (0,0,0)^T$. More precisely the p.d.f. of $\tau(\omega, z)$ is

$$p(\tau) = \frac{\tau^{1/2}}{\sqrt{2\pi}(\sigma^2 z)^{3/2}} \exp\left(-\frac{\tau}{2\sigma^2 z}\right) \mathbf{1}_{[0,\infty)}(\tau).$$

5.2.2. *Realistic models.* The white noise model gives an analytical formula for the p.d.f. of the DGD, which in turns allows us to compute exactly the probability that the DGD exceeds a given threshold value. However, it has been pointed out that the p.d.f. tail of the DGD does not fit with the Maxwellian distribution in realistic configurations [1]. Various numerical and experimental PMD generation techniques involve the concatenation of birefringent elements with piecewise constant vectors $\Omega$ [18]. Equation (5.1) can be solved over each segment, and continuity conditions on the segments junctions give a discrete model for the PMD vector $\hat{\mathbf{r}}$. The total PMD vector after the $(n+1)$st section can then be obtained from the concatenation equation [14]

$$\hat{\mathbf{r}}_{n+1} = R_{n+1}\hat{\mathbf{r}}_n + \sigma\Omega_{n+1}, \tag{5.5}$$

where $\sigma$ is the DGD per section, $\Omega_n = \Omega(\theta_n)$ with

$$\Omega(\theta) = (\cos(\theta), \sin(\theta), 0)^T.$$

$R_n$ is a matrix corresponding to a rotation through an angle $\phi_n$ about the axis $\Omega_n$. Explicitly $R_n = R(\theta_n, \phi_n)$ with

$$R(\theta, \phi) = \begin{pmatrix} \cos^2(\theta) + \sin^2(\theta)\cos(\phi) & \sin(\theta)\cos(\theta)(1-\cos(\phi)) & \sin(\theta)\sin(\phi) \\ \sin(\theta)\cos(\theta)(1-\cos(\phi)) & \sin^2(\theta) + \cos^2(\theta)\cos(\phi) & -\cos(\theta)\sin(\phi) \\ -\sin(\theta)\sin(\phi) & \cos(\theta)\sin(\phi) & \cos(\phi) \end{pmatrix}.$$

From the probabilistic point of view, the angles $\phi_n$ are i.i.d. random variables uniformly distributed in $(0, 2\pi)$. The angles $\theta_n$ are i.i.d. random variables such that $\cos(\theta_n)$ are uniformly distributed in $(-1, 1)$ [2]. Accordingly, $(\hat{\mathbf{r}}_n)_{n \in \mathbb{N}}$ is a Markov chain. Let us assume that the fiber is modeled as the concatenation of $n$ segments and that the outage event is of the form $|\hat{\mathbf{r}}_n| > a$ for some fixed threshold value $a$. In the case where $a$ is much larger than the expected value of the final DGD $|\hat{\mathbf{r}}_n|$, the outage probability is very small, and this is the quantity that we want to estimate.

### 5.3. *Estimations of outage probabilities.*

5.3.1. *Importance sampling.* In [2] IS is used to accurately calculate outage probabilities due to PMD. As discussed in the Introduction, the key difficulty in applying IS is to properly choose the twisted distribution for the driving process $(\theta_p, \phi_p)_{1 \leq p \leq n}$. The papers [2, 11, 17] present different twisted distributions and the physical explanations why such distributions are likely to produce large DGD's. As a result the authors obtain with $10^5$ simulations good approximations of the p.d.f. tail even for probabilities of the order of $10^{-12}$. The main reported physical result is that the probability tail is significantly smaller than the Maxwellian tail predicted by the white noise model.

5.3.2. *Interacting particle systems.* In this subsection we apply our IPS method and compare the results with those obtained by MC and IS. To get a reliable estimate of the outage probability of the event, it is necessary to generate realizations producing large DGD's. The main advantage of the IPS approach is that it proposes a "blink" method that does not require any physical insight. Such a method could thus be generalized to more complicated situations. Here the Markov process is the PMD vector $(\hat{\mathbf{r}}_n)_{n \in \mathbb{N}}$ at the output of the $n$th fiber section. The state space is $\mathbb{R}^3$, the initial PMD vector is $\hat{\mathbf{r}}_0 = (0, 0, 0)^T$, the Markov transitions are described by (5.5) and the energy-like function is $V(\hat{\mathbf{r}}) = |\hat{\mathbf{r}}|$. We estimate the p.d.f. $p(a)$ of $|\hat{\mathbf{r}}_n|$ by implementing the IPS with the two weight functions

$$(5.6) \qquad G_p^\beta(\hat{\mathbf{r}}) = \exp(\beta|\hat{\mathbf{r}}_p|)$$

parameterized by $\beta \geq 0$, and

(5.7) $$G_p^\alpha(\hat{\mathbf{r}}) = \exp[\alpha(|\hat{\mathbf{r}}_p| - |\hat{\mathbf{r}}_{p-1}|)]$$

parameterized by $\alpha \geq 0$. We have implemented Algorithms 1 and 2 as described in Section 3. Before presenting and discussing the results, we would like to underline that we have chosen to perform a selection step at the output of each segment, because the number of segments is not very large. If the number of segments were very large, it should be better to perform a selection step every two or three or $n_0$ segments.

In Figure 3(a) we plot the estimation of the DGD p.d.f. obtained by the IPS method with the weight function $G_n^\beta$ defined by (5.6). The fiber consists in the concatenation of $n = 15$ segments. The DGD per section is $\sigma = 0.5$. We use a set of $N = 2 \times 10^4$ interacting particles. This result can be compared with the one obtained in [2], which shows excellent agreement. The difference is that our procedure is fully adaptative, it does not require any guess of the user, and it does not require to twist the input probability density. The variance $p_2^2$ of the estimator of the DGD p.d.f. is plotted in Figure 3(b). This figure is actually used to determine the best estimator of the DGD p.d.f. by the procedure described at the end of Section 3.

In Figure 4(a) we plot the estimation of the DGD p.d.f. obtained by the IPS method with the weight function $G_p^\alpha$ defined by (5.7). It turns out that the estimated variance of the estimator is smaller with the weight function $G_p^\alpha$ than with the weight function $G_p^\beta$ [cf. Figures 4(b) and 3(b)]. This observation confirms the theoretical predictions obtained with the Gaussian random walk.

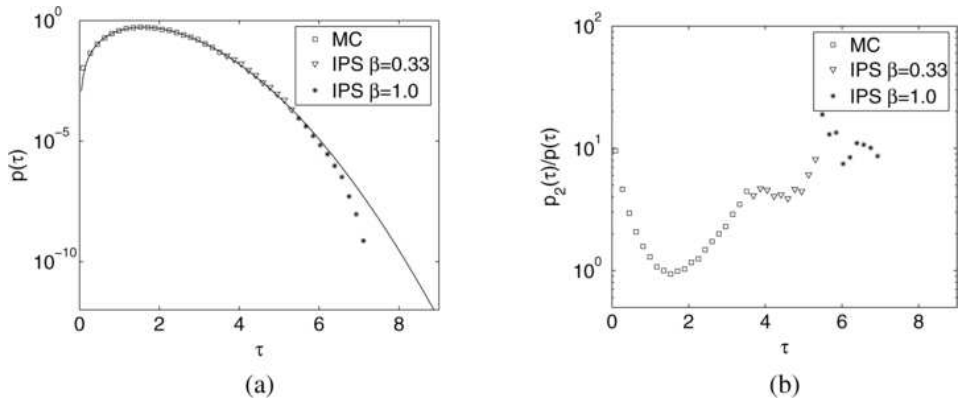

Fig. 3. (a) *Segments of the DGD p.d.f. obtained by the usual MC technique (squares) and by the IPS with the weight function $G_n^\beta$ with $\beta = 0.33$ (triangles) and $\beta = 1$ (stars). The solid line stands for the Maxwellian distribution obtained with the white noise model. The Maxwellian distribution fails to describe accurately the p.d.f. tail.* (b): *Ratios $p_2/p$. The quantity $p_2$ is the standard deviation of the estimator of the DGD p.d.f. In the IPS cases, the standard deviations are estimated via the formula* (3.17).

GENEALOGICAL PARTICLE ANALYSIS OF RARE EVENTS 39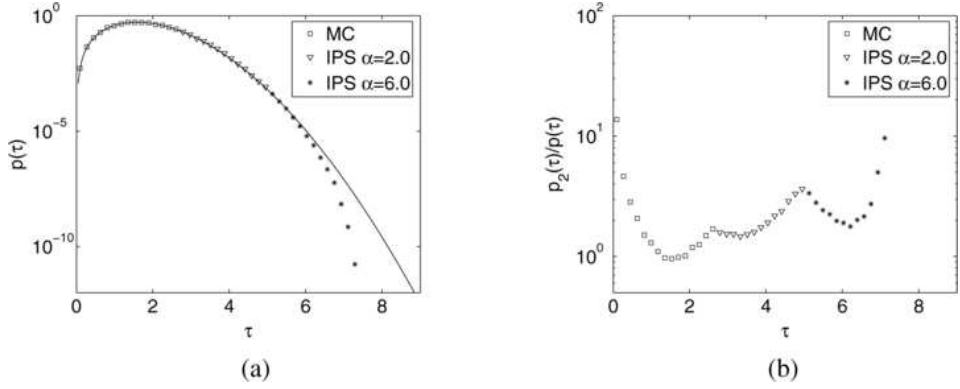

FIG. 4. (a) *Segments of the DGD p.d.f. obtained by the usual MC technique (squares) and by the IPS with the weight function $G_p^\alpha$ with $\alpha = 2.0$ (triangles) and $\alpha = 6.0$ (stars). The solid line stands for the Maxwellian distribution obtained with the white noise model.* (b) *Ratios $p_2/p$. The quantity $p_2$ is the standard deviation of the estimator of the DGD p.d.f. In the IPS cases, the standard deviations are estimated via* (3.18).

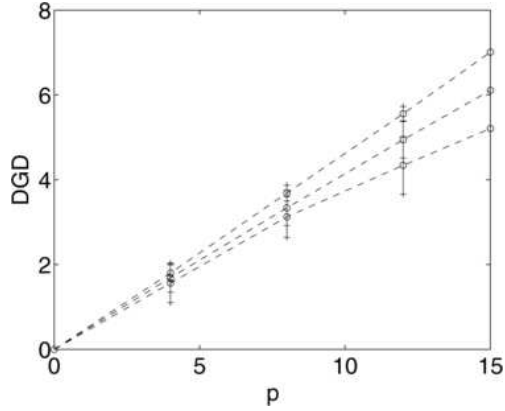

FIG. 5. *Conditional expectations $D_{a,a+\delta a}^1(p,n)$ of the intermediate DGD at $p = 4$, 8, 12, given that the final DGD lies in the interval $(a, a + \delta a)$ with $n = 15$, $\delta a = 0.18$, and (from top to bottom) $a = 7$, $a = 6.1$, $a = 5.2$. The error bars are obtained from the estimations of the conditional variances.*

The IPS approach is also powerful to compute conditional probabilities or expectations given the occurrence of some rare event. For instance, we can be interested in the moments of the intermediate DGDs given that the final DGD lies in the rare set $(a, a + \delta a)$:

$$D_{a,a+\delta a}^q(p,n) = \mathbb{E}[|\hat{\mathbf{r}}_p|^q | |\hat{\mathbf{r}}_n| \in [a, a + \delta a)].$$

This information gives us the typical behaviors of the PMD vectors along the fiber that give rise to a large final DGD. We use the estimator (2.20) based on the IPS with the weight function (5.7). As shown by Figure 5, the typical conditional trajectory of the DGD is close to a linear increase with a constant rate given by the ratio of the final DGD over the length of the fiber. The conditional variances are found to be small, which shows that fluctuations are relatively small around this average behavior.

Laboratoire de Mathématiques J. A. Dieudonné
Université de Nice
Parc Valrose
06108 Nice Cedex 2
France
e-mail: delmoral@math.unice.fr

Laboratoire de Probabilités
et Modèles Aléatoires
and Laboratoire Jacques-Louis Lions
Université Paris 7
2 Place Jussieu
75251 Paris Cedex 5
France
e-mail: garnier@math.jussieu.fr